\documentclass[12pt,reqno]{amsart}
\usepackage{amsthm}
\usepackage{amssymb}
\usepackage{amsmath}
\usepackage{mathrsfs}
\usepackage{amsfonts}
\usepackage{amssymb,amsmath}
\usepackage{graphicx}  
 \usepackage[colorlinks, linkcolor=black, citecolor=blue,   hypertexnames=false]{hyperref}
 \usepackage[numbers,sort&compress]{natbib}
 \usepackage{epsfig}
 \usepackage{subfigure}
\newtheorem{Theorem}{Theorem}[section]

\newtheorem{Lemma}{Lemma}[section]
\newtheorem{Corollary}{Corollary}[section]

\newtheorem{Definition}{Definition}[section]
\numberwithin{equation}{section}

\numberwithin{equation}{section}

\def \no#1#2#3 {{\bf #1} (#3), #2.}
\def \eds#1#2#3 {#1, #2, #3.}

\title[Asymptotic Stability of Solitons to 1D NLS]
{\bf{
Asymptotic stability of solitons to 1D Nonlinear Schr\"odinger Equations in the subcritical case}}

\author[Z. Li ]
{Ze Li}

\address{Ze Li
\newline\indent
School of Mathematics and Statistics, Ningbo University
\newline\indent
Ningbo, 315000, Zhejiang, P.R. China
}
\email{rikudosennin@163.com}

\keywords{nonlinear Schr\"odinger equations, solitons, asymptotic stability, weak nonlinearity}

\begin{document}

\begin{abstract}
In this paper, we prove the asymptotic stability of solitary waves to 1D nonlinear Schr\"odinger equations in the subcritical case with symmetry and spectrum assumptions. One of the main ideas is to use the vector fields method developed by Cuccagna, Georgiev, Visciglia \cite{CGV} to overcome the weak decay with respect to $t$ of the linearized equation caused by the one dimension setting and the weak nonlinearity caused by the subcritical growth of the nonlinearity term. Meanwhile, we apply the polynomial growth of the high Sobolev norms of solutions to 1D Schr\"odinger equations obtained by Staffilani \cite{Sta} to control the high moments of the solutions  emerging from the vector fields method.
\end{abstract}

\maketitle

\section{Introduction}
In this paper, we consider the following nonlinear Schr\"odinger equation (NLS),
\begin{align}\label{1}
\left\{ \begin{array}{l}
 i{u_t} +\Delta u = F(\left| u \right|^2)u \\
 u(0,x) = {u_0}(x).
 \end{array} \right.
\end{align}
where $u:\Bbb R\times \Bbb R\to \Bbb C$.
The nonlinear Schr\"odinger equation is a classical field equation whose principal applications are to the propagation of light in nonlinear optical fibers and planar waveguides (e.g. \cite{MB}), and to Bose-Einstein condensates confined to highly anisotropic cigar-shaped traps in the mean-field regime (e.g. \cite{PS}).

NLS has a family of localized solutions among which the so-called solitons are the best understood.  We are interested in the dynamics of NLS around solitons. The orbital stability of solitons and multi-solitons, was considered by many authors for various models, for instance Benjamin \cite{Be}, Grillakis, Shatah, Strauss \cite{GSS1,GSS2}, Weinstein \cite{W2}, Martel, Merle, Tsai \cite{MMT}. The other main problem on solitons is the asymptotic stability
which states that any solution of NLS initiated near the family of solitons decomposes into a moving soliton and a radiation part with an asymptotically vanishing remainder as $t\to\infty$. This is best known for completely integrable equations for instance the one dimensional cubic NLS by using the inverse scattering method. For general nonlinearities, the first asymptotic stability result was obtained by Soffer, Weinstein \cite{SW} in context of the equation
\begin{align}\label{2}
iu_t+\Delta u +V(x)u=F(|u|^2)u.
\end{align}
There have been a lot of works in the study of asymptotic stability of solitons for (\ref{2}) especially for attractive potentials, for instance \cite{CM1,GNT,T,KZ,KM1}.

The asymptotic stability of solitons for (\ref{1}) started from Buslaev, Perelman \cite{BP} where the one dimension NLS was considered and further refined in \cite{BP2,BP3}. Their work was extended to high dimensions in Cuccagna \cite{C}. Later Perelman \cite{P}, Rodnianski, Schlag, Soffer \cite{RSS} proved the asymptotic stability for multi-solitons in high dimensions. When blow up solution exists such as the super-mass critical and mass critical NLS equations, instead of asymptotic stability, stable and center stable manifolds are introduced to describe the dynamics of NLS near solitons, see Krieger, Schlag \cite{KS,KS2}, Beceanu \cite{B1,B2,B3}, Cuccagna \cite{C2}, Schlag \cite{S}.

The idea in this paper to prove asymptotic stability of solitons for (\ref{1}) is to use the vector field method to obtain more decay of the solution to the linearized equation with respect to $t$. This appeared originally in Klainerman \cite{K1,K2} where the small data global well-posedness of quasilinear wave equations was solved. One of the key idea of the vector field method is to gain the decay with respect to time by paying more decay and regularity with respect to the spatial variables.  The same idea works for Schr\"odinger equations see McKean, Shatah \cite{MS} for NLS and  Cuccagna, Georgiev, Visciglia \cite{CGV} for NLS with a potential.

In order to state our theorem, we give the definition of solitary waves.
\begin{Definition}
We call the periodic localized solution to (\ref{1}),
$$w(x;\sigma)=exp(-i\beta+i\frac{1}{2}vx)\varphi(x-b;\alpha),$$
solitary wave, if $\varphi(x;\alpha)$ is a radial positive function and the time dependent parameters  $\sigma(t)\triangleq(\beta(t),\omega(t),b(t),v(t))$ satisfy
\begin{align}
&\Delta_x\varphi=\frac{\alpha^2}{4}\varphi+F(\varphi^2)\varphi\nonumber\\
&\beta'=\omega, \omega'=0, b'=v, v'=0, \omega=\frac{1}{4}(v^2-\alpha^2).\label{70}
\end{align}
\end{Definition}
It is known that when $F(x)=|x|^{\frac{p-1}{2}}$, $p$ is mass-subcritical, $\varphi$ exists and decays exponentially at infinity. (See Berestycki, Lions \cite{BL} where a larger class of nonlinearity was considered.)

\subsection{Linearized operator}
As in \cite{BP}, the linearization operator of (1.1) around the solitary wave $w(x,t;\sigma)$ is given by
\begin{align*}
i\partial_t{\chi}=-\Delta\chi+F(|w|^2)\chi+F'(|w|^2)w(w\chi+w\overline {{\chi}}).
\end{align*}
If we denote
\begin{align*}
\chi(x,t)=e^{i\Phi}f(y,t), \mbox{  }\Phi=-\beta(t)+\frac{1}{2}vx, \mbox{   }y=x-b,
\end{align*}
then the function $f$ satisfies
\begin{align*}i\partial_t{f}=L(\alpha)f,
\end{align*}
where
\begin{align}\label{app}
L(\alpha)f=-\Delta f+\alpha^2f/4+F(\varphi^2)f +F'(\varphi^2)\varphi^2(f+\overline f), \varphi=\varphi(y,\alpha).
\end{align}
Consider the complexification of (\ref{app}):
{\small\begin{align*}
&i\partial_t{ \mathbf{f}}=H(\alpha){ \bf f}, {\bf f}=(f,\overline{{f}})^t,\\
&H(\alpha)=H_0(\alpha)+V(\alpha), H_0(\alpha)=(-\Delta_y+\alpha^2/4)\theta_3,\\
&V(\alpha)=[F(\varphi^2)+F'(\varphi^2)\varphi^2]\theta_3+iF'(\varphi^2)\varphi^2\theta_2,
\end{align*}}
where $\theta_2$ and $\theta_3$ are the Pauli matrices:
{\small\begin{align*}
{\theta _2} = \left( \begin{array}{l}
 0 \\
 i \\
 \end{array} \right.\left. \begin{array}{l}
  - i \\
 0 \\
 \end{array} \right),  \mbox{  }{\theta _3} = \left( \begin{array}{l}
 1 \\
 0 \\
 \end{array} \right.\left. \begin{array}{l}
 0 \\
  - 1 \\
 \end{array} \right).
\end{align*}}

There are four known generalized eigenfunctions denoted by $\{\xi_1,\xi_2,\xi_3,\xi_4\}$ for the zero eigenvalue to $H(\alpha)$.  Among them the two radial functions are
{\small\begin{align*}
{\xi _1} = \left( \begin{array}{l}
 {v_1} \\
 {{\bar v}_1} \\
 \end{array} \right),{\xi _2} = \left( \begin{array}{l}
 {v_2} \\
 {{\bar v}_2} \\
 \end{array} \right),{v_1} =  - i\varphi (y;\alpha ),{v_2} =  - \frac{2}{\alpha }{\varphi _\alpha }(y;\alpha ).
\end{align*}}
Moreover, $H(\alpha)\xi_1=0$, $H(\alpha)\xi_2=i\xi_1$, $\left\langle {{\xi _1},{\xi _2}} \right\rangle  = 0$.

\subsection{Main results}
We have the following assumptions about the nonlinearity $F$.\\
{\bf Assumption A}:\\
\noindent$(i)$ There exists a solitary wave solution to (\ref{1}) and it is of exponential decay.\\
\noindent$(ii)$ There exist $m,n\ge 4$ such that for $k\in \{0,1,2,3,4\}$,
{\footnotesize\begin{align*}
&\left| {\frac{{{d^k}}}{{d{x^k}}}F(x)} \right| \le {\left| x \right|^{\frac{{m - 1}}{2} - k}},\mbox{  }{\rm{if}} \mbox{  }\left| x \right| \le 1, \\
&\left| {\frac{{{d^k}}}{{d{x^k}}}F(x)} \right| \le {\left| x \right|^{\frac{{n - 1}}{2} - k}},\mbox{  }{\rm{if}} \mbox{  }\left| x \right| \ge 1.
\end{align*}}
\noindent$(iii)$ The linearized operator $H(\alpha)$ has zero as its only eigenvalue with generalized eigenfunction space spanned by $\{\xi_1, \xi_2, \xi_3, \xi_4\}$, no resonance and no imbedded eigenvalues in its continuous spectrum. \\
\noindent$(iv)$ (\ref{1}) is globally well-posed in $H^1$ and the solution $u(t)$ satisfies $\|u(t)\|_{H^1}\le C(\|u_0\|_{H^1})$.

\noindent{\bf Remark 1.1} If we restrict ourself to the linear combination of power functions, roughly speaking, $(ii)$ means the lowest degree of $F$ is $\frac{3}{2}$, which improves the result in \cite{BP} where the degree of $F$ is assumed to be at least four.\\
{\bf Remark 1.2} For $n<5$, it is proved in Weinstein \cite{W1,W2} that the generalized eigenfunction space to zero of $H(\alpha)$ is four dimensional. For cubic nonlinearity in three dimensions, Costin, Huang, Schlag \cite{CHS} proved that $H(\alpha)$ has no imbedded eigenvalues. However, whether there exists imbedded eigenvalue or eigenvalues in the gap between zero and the continuous spectrum in the one dimension case for subcritical power nonlinearity is unknown.\\
{\bf Remark 1.3} Although $(iii)$ is widely assumed in the papers studying asymptotic stability, the rationality of $(iii)$ is unknown when $F(x)=x^p$ and $0<p<2$. However there are some possible remedies to deal with the case when more than one discrete spectrum occurs, such as the so called `` Fermi Golden Rule" hypothesis (FGR) introduced by Sigal \cite{S1}. See also Cuccagna \cite{C2}, Soffer, Weinstein \cite{SW1,SW2} for applications of normal forms and Fermi Golden Rule in the study of dynamics near solitons.\\
{\bf Remark 1.4} If we assume $n<5$, then $(iv)$ is naturally satisfied (e.g.\cite{Taz}). Meanwhile, $(iv)$ may hold for combined nonlinearities.

Our main theorem is as follows.
Let $\|f\|_{\Sigma}=\|f\|_{H^2_x}+\|xf\|_{L^2_x}+\|x^2f\|_{L^2_x}$.
\begin{Theorem}
Let $F$ satisfy Assumption A.
Assume that $w(x;\sigma_0(t))$ is a solitary wave solution to (\ref{1}) with $\sigma_0(t)=(\beta_0(t),\omega_0(t),0,0)$ and satisfies
$\frac{d}{d\alpha}\|\varphi\|_2^2\neq 0$ at $\alpha_0$, where $-\frac{1}{4}\alpha_0^2=\omega_0$. If the radial initial data $u_0$ satisfies
\begin{align}\label{p0kinbgh}
\| u_0-w(x;\sigma_0(0)) \|_{\Sigma} \ll1,
\end{align}
then there exists a modulated solitary wave $w(x;\sigma_+(t))$ such that as $t\to\infty$ $u(x,t)$ decomposes into
$$u(x,t)=w(x;\sigma_+(t))+\chi(t),$$
where $\|\chi(t)\|_{\infty}\le Ct^{-s+1}$ with $s=\frac{7}{4}^+$.
\end{Theorem}

\subsection{Main ideas}

The key ingredient of our proof is to establish the decay estimate with some proper vector field operator $|J_V(t)|^s$,
{\footnotesize\begin{align}\label{4}
\|\chi\|_{\infty}\le Ct^{-s} \||J_V(t)|^s\chi\|_{L^{\infty}_tL^{2}_x}.
\end{align}}
Let $\mathcal{H}(\alpha)=P_{c}(H(\alpha))H(\alpha)$ be the projection of $H(\alpha)$ onto its continuous spectrum space.  We will choose $|J_V(t)|^s$ to be $U(t+h)t^s\mathcal{H}(\alpha(t))^{\frac{s}{2}}U(-t-h)$, see Section 3.
The importance of (\ref{4}) is that it fills the gap between the weak decay resulting from dispersive effects in the one dimension and the desired decay to put the quadratic remainder in $L^1_t$.

There are two main difficulties in the establishment of (\ref{4}). The first is how to give a proper definition of $\mathcal{H}(\alpha)^{\frac{s}{2}}$.
We notice that it is not trivial since $\mathcal{H}(\alpha)$ is not selfadjoint.
The second difficulty is how to prove (\ref{4}). A similar estimate for Schr\"odinger operators is given in \cite{CGV} by the inverse scattering theory. This may not work for our $\mathcal{H}(\alpha)$ since it is not selfadjoint and no inverse scattering is available for it.
To overcome aforementioned difficulties, we rely heavily on the functional calculus techniques. Specifically speaking, the fractional power of $\mathcal{H}(\alpha)$ is defined via Dunford-Schwartz type integral with a properly chosen integral contour (see formula (\ref{20}) ).  We emphasize again the fractional power is defined for $\mathcal{H}(\alpha)$ rather than ${H}(\alpha)$ itself.  Thus we need to separate the discrete spectrum part from the solution and study the modulation equations. To get (\ref{4}), we translate the problem to the corresponding estimate of the inverse of $\mathcal{H}(\alpha)^{\frac{s}{2}}$ which finally reduces to estimates of resolvent $(H(\alpha)-\lambda)^{-1}$. When $\lambda$ is sufficiently large, we apply a perturbation technique to obtain the decay of the resolvent $(\lambda-H(\alpha))^{-1}$.

It is important that ({\ref{4}) remains valid only for $\chi$ satisfying $U(-t-h)\chi\in P_c(H(\alpha))$. Hence orthogonal conditions should be added to $\chi$. As soon as we establish (\ref{4}), the main theorem follows by the study of modulation equation and Strichartz estimates.

The paper is organized as follows. In Section 2, we define the fractional power of the linearized operator. In Section 3, we define the vector field operator $|J_V|^s$ and introduce corresponding orthogonal conditions.  In Section 4, we set up the bootstrap argument and study the modulation equation.  In Section 5, we derive the equation of $|J_V|^su$. In Section 6, we close the bootstrap by nonlinear estimates and prove Theorem 1.1.

{\bf Notations} $A\lesssim B$ means there exists some positive constant $C$ such that $A\le CB$. $A\sim B$ means $A\lesssim B$ and $B\lesssim A$. Without confusion $\|\cdot\|_{L^p_x}$ will always be denoted as $\|\cdot\|_{p}$. The spacetime norm $(\int^{s_2}_{s_1}(\int_{\Bbb R}|f(t,x)|^pdx)^\frac{r}{p}dt)^{\frac{1}{r}}$ is denoted by $\|f\|_{L^r{[s_1,s_2]}L^p_x}$.

\section{Fractional power of $H(\alpha)$ in continuous spectrum space}
In this section we give the definition of the fractional power of $H(\alpha)$. And some distorted Sobolev embedding and equivalence lemmas are proved.
Given $\alpha>0$, recall the linearized operator $H(\alpha)$. We will work on $L^2_{rad}$, and the domain of $H(\alpha)$ is taken as $D(H(\alpha))=H^2_{rad}$.
 Define $\tau=\frac{1}{4}\alpha^2$ and
$${H_0} = \left( \begin{array}{l}
  - \Delta  + \tau  \\
  \\
 \end{array} \right.\left. \begin{array}{l}
  \\
 \Delta  - \tau  \\
 \end{array} \right).
$$
With these notations, the potential term in $H(\alpha)=H_0+V$ is of the form
$$\left( \begin{array}{l}
 {V_1} \\
 {V_2} \\
 \end{array} \right.\left. \begin{array}{l}
  - {V_2} \\
  - {V_1} \\
 \end{array} \right).
 $$

The explicit formula for the Green function of $H_0$ and Young's inequality yield
\begin{Lemma}\label{45}
For $1\le p\le\infty$, $\lambda\in {\Bbb C}\backslash[(-\infty, -\tau]\bigcup[\tau, \infty)]$,
\begin{align*}
&{\left\| {{{({H_0} - \lambda )}^{ - 1}}} \right\|_{p \to p}} \le {c_1}\big[{({\Re}\sqrt { \tau-\lambda } )^{ - 2}} + {({\Re}\sqrt { \lambda+ \tau } )^{ - 2}}\big];
\end{align*}
and for $2\le p\le \infty$,
\begin{align*}
{\left\| {{{({H_0} - \lambda )}^{ - 1}}} \right\|_{2 \to p}} \le {c_1}\big[{({\Re}\sqrt {\tau-\lambda} )^{ - \frac{3}{2} - \frac{1}{p}}} + {({\Re}\sqrt { \lambda+\tau } )^{ - \frac{3}{2} - \frac{1}{p}}}\big].
\end{align*}
\end{Lemma}

The following is the high frequency estimate.
\begin{Lemma}\label{41}(High frequency resolvent estimate)
Let $\lambda\in {\Bbb C}\backslash\big[(-\infty, -\tau]\bigcup[\tau, \infty)\big]$.
Suppose that $\lambda$ is chosen such that $({\Re}\sqrt{\tau-\lambda})^2\ge4c_1\|V\|_{\infty}$ and $({\Re}\sqrt{\lambda+\tau})^2\ge4c_1\|V\|_{\infty}$, $1\le p\le \infty$, then
\begin{align}\label{hulijin}
{\|(H(\alpha) - \lambda )^{ - 1}}\|{_{p \to p}} \le c\big[{({\Re}\sqrt {\tau-\lambda} )^{ - 2}} + {({\Re}\sqrt { \lambda+\tau } )^{ - 2}}\big],
\end{align}
where $c$ is independent of $V$.
\end{Lemma}
\begin{proof}
By Lemma \ref{45},
\begin{align*}
{\left\| {{{({H_0} - \lambda )}^{ - 1}}V} \right\|_{p \to p}} &\le {\left\| {{{({H_0} - \lambda )}^{ - 1}}} \right\|_{p \to p}}{\left\| V \right\|_\infty }\\
& \le {c_1}[{({\Re}\sqrt {\tau-\lambda} )^{ - 2}}+ {({\Re}\sqrt { \lambda+\tau } )^{ - 2}}]{\left\| V \right\|_\infty }.
\end{align*}
Hence if
${({\Re}\sqrt {\tau-\lambda} )^2} \ge 4{c_1}\|V\|_{\infty }, \mbox{  }{({\Re}\sqrt { \lambda+\tau } )^2} \ge 4{c_1}\|V\|_{\infty }
$, we have
$${\left\| {{{(H(\alpha) - \lambda )}^{ - 1}}} \right\|_{p \to p}} \le 2{\left\| {{{({H_0} - \lambda )}^{ - 1}}} \right\|_{p \to p}} \le c[{({\Re}\sqrt {\tau-\lambda} )^{ - 2}} + {({\Re}\sqrt { \lambda+\tau } )^{ - 2}}].$$
Thus ${\left\| {{{({H_0} - \lambda )}^{ - 1}}V} \right\|_{p \to p}} < \frac{1}{2}$. Then by Neumann series,
$${(H(\alpha) - \lambda )^{ - 1}} = {(I + {({H_0} - \lambda )^{ - 1}}V)^{ - 1}}{({H_0} - \lambda )^{ - 1}},$$
and (\ref{hulijin}) follows.
\end{proof}

{\bf Denote $P_c({H}(\alpha))$  the projection onto continuous spectrum space of  ${H}(\alpha)$ in $L^2_{rad}$. And let $\mathcal{H}(\alpha)=H(\alpha)P_c({H}(\alpha))$.}
\noindent For $0<s<2$, $y\in H^2$, define
\begin{align}\label{20}
{\mathcal{H}(\alpha)^{\frac{s}{2}}}y = \frac{1}{2\pi i}\int_\Gamma  {{\lambda ^{\frac{s}{2} - 1}}{{(\lambda  - \mathcal{H}(\alpha))}^{ - 1}}\mathcal{H}(\alpha)yd\lambda }, \mbox{ }0<s<2, y\in H^2,
\end{align}
where
$$\Gamma  = {\gamma _{1, + }} \cup {\gamma _{2, + }} \cup {\gamma _{1, - }} \cup {\gamma _{2, - }},\mbox{  }{\gamma _{1, \pm }}(t) = t + a \pm i\varepsilon t, \mbox{  }{\gamma _{2, \pm }}(t) =  - t - a \pm i\varepsilon ,t \ge 0.
$$
{\bf We emphasize ${\mathcal{H}}(\alpha)^{\frac{s}{2}}$ is not ${H}(\alpha)^{\frac{s}{2}}$ in general}.

See Figure 1 for the shape of $\Gamma$ and we make the convention that the two connected branches of $\Gamma$ are anticlockwise oriented.
Here $a>0$ is some appropriate constant close to $\tau$ excluding the discrete spectrum of $H(\alpha)$ in the interior of $\Gamma$, and $\varepsilon>0$ is  to ensure $ {\gamma _{1, \pm }}(t)$ and $ {\gamma _{2, \pm }}(t)$ lie in some single-value branch of $z^{\frac{s}{2}}$.

The integration in (\ref{20}) is well-defined for each $y\in H^2$ because of Lemma \ref{41}. In fact, the analyticity of $(H(\alpha)-\lambda)^{-1}$ in the resolvent set and the fact $\Gamma$ is strictly away from the discrete spectrum $\{0\}$ show (\ref{20}) is well defined for $\lambda$ contained in a ball. Meanwhile, Lemma \ref{41} shows the integrand in (\ref{20}) decays fast as $\lambda$ goes to infinity along $\Gamma$. Thus (\ref{20}) is well defined for $y\in H^2$.

For $\beta<0$, define
\begin{align}\label{Po89nb}
{[[\mathcal{H}(\alpha)]]^\beta }y= \frac{1}{2\pi i}\int_\Gamma  {{\lambda ^\beta }{{(\lambda  - \mathcal{H}(\alpha))}^{ - 1}}y} d\lambda, \mbox{ }\beta<0.
\end{align}
For each $y\in L^2$, the integration is well-defined due to Lemma \ref{41}. {\bf The notation $[[\cdot]]$ in (\ref{Po89nb}) doesnot mean brackets, we use it to distinguish $[[\mathcal{H}(\alpha)]]^\beta$ from  $\mathcal{H}(\alpha)^s$ in (\ref{20}).}

Similarly, we define $\mathcal{H}_0^{\frac{s}{2}}$ with $0<s<2$ by
\begin{align}\label{109b}
{\mathcal{H}_0^{\frac{s}{2}}}y = \frac{1}{2\pi i}\int_\Gamma  {{\lambda ^{\frac{s}{2} - 1}}{{(\lambda  - {H}_0)}^{ - 1}}{H}_0yd\lambda }, y\in H^2,
\end{align}
and $[[\mathcal{H}_0]]^{\beta}$ with $\beta<0$ to be
\begin{align}\label{10000b}
{[[\mathcal{H}_0]]^\beta }y= \frac{1}{2\pi i}\int_\Gamma  {{\lambda ^\beta }{{(\lambda  -  {H}_0)}^{ - 1}}y} d\lambda,  y\in L^2.
\end{align}
Notice that since ${H}_0$ has only continuous spectrum, there is no needs to use $P_c(H_0)$ in the definition.

\begin{Lemma}\label{xinhui}
For $0<s<2$, there holds
\begin{align}
{\mathcal{H}(\alpha)}^{\frac{s}{2}}\xi_i(\alpha)=0, i=1,2.
\end{align}
\end{Lemma}
\begin{proof}
This follows by $\mathcal{H}(\alpha)=P_{c}(H(\alpha))$ and $\{\xi_i(\alpha)\}_{i=1,2}\subset H^2$ are generalized eigenfunctions.
\end{proof}

\begin{figure}[htbp]
\centering
\subfigure[curve $\Gamma$]{
\begin{minipage}[t]{0.3\linewidth}
\centering
\includegraphics[width=4.3cm,height= 4.3cm]{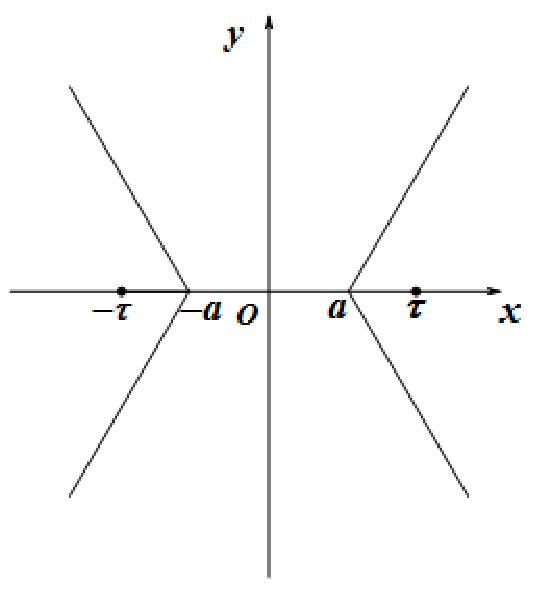}
\end{minipage}
}
\subfigure[curve $\mathcal{C}_n$]{
\begin{minipage}[t]{0.3\linewidth}
\centering
\includegraphics[width=5.2cm,height= 4cm]{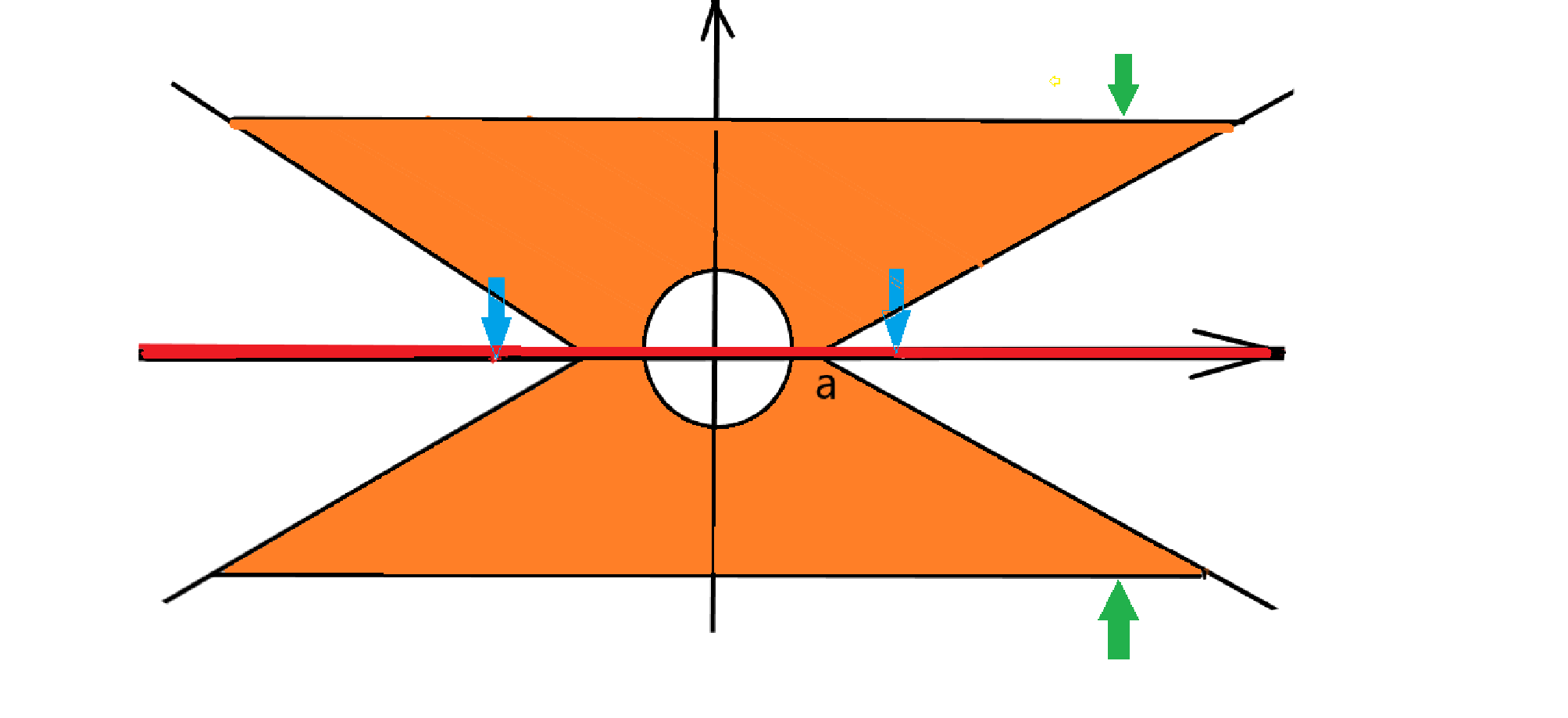}
\end{minipage}
}

\quad
\subfigure[curve $\Gamma$ \mbox{  }and\mbox{  }$\Gamma'$ ]{
\begin{minipage}[t]{0.3\linewidth}
\centering
\includegraphics[width=4.7cm,height= 4.7cm]{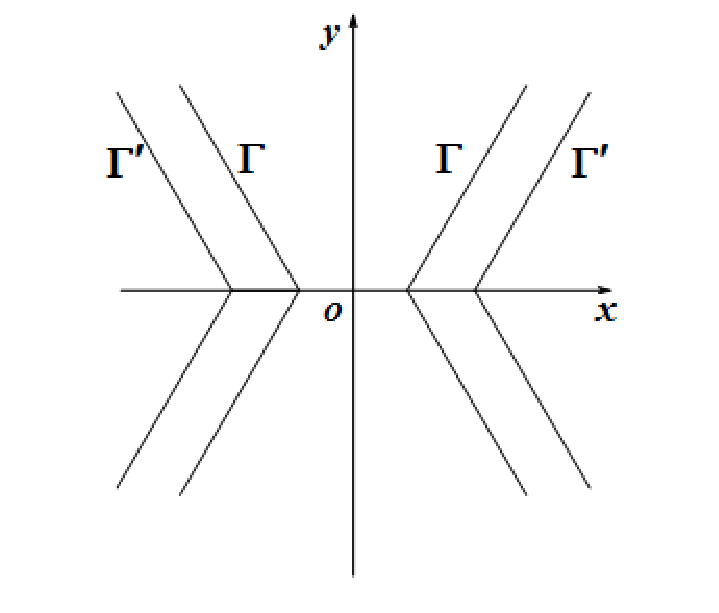}
\end{minipage}
}
\subfigure[{curve $\Gamma_{1,n}$, $\Gamma_{1,n}'$, and $\gamma_{1,n}'$ }.]{
\begin{minipage}[t]{0.3\linewidth}
\centering
\includegraphics[width=5cm,height= 5cm]{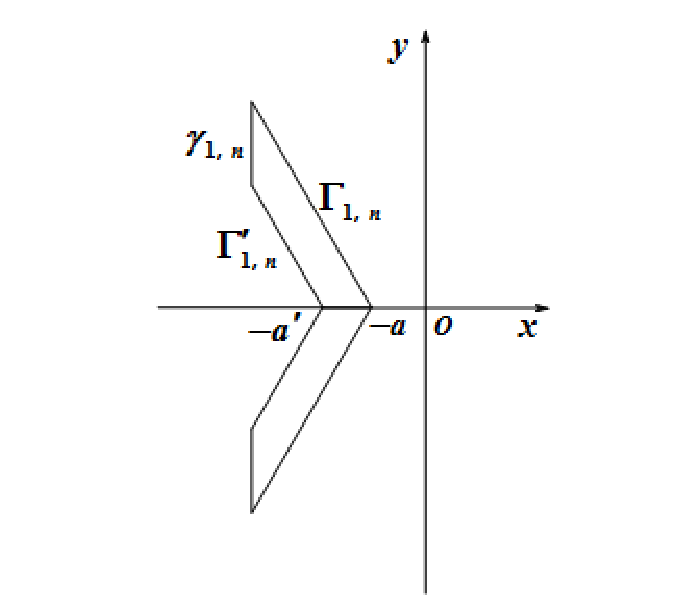}
\end{minipage}
}
\end{figure}

\noindent{\bf Remark 2.1} Let us explain the logic behind the integral contours in above pictures. In order to make $f(z)=z^{\theta}$ be  single valued analytic functions, we remove the negative $y$-axis. That is why we choose integral contours like (a), (c), (d).
However, in the special case $f(z)=z^{-1}$ which is analytic outside 0, integral contours like (b) are reasonable.

\begin{Lemma}\label{51a}
Denote ${\rm Range}( P_{c}(H(\alpha)))$ the continuous spectrum space of $H(\alpha)$ in $L^2_{rad}$. We have
\begin{align*}
\mathcal{H}(\alpha)^{\frac{s}{2}}g=([[\mathcal{H}(\alpha)]]^{-\frac{s}{2}})^{-1}g, \mbox{ } \forall g \in H^2\cap {\rm Range}( P_{c}(H(\alpha))).
\end{align*}
\end{Lemma}
\begin{proof}
{\bf{Step 1}}. In this step, we prove $[[\mathcal{H}(\alpha)]]^{-1}$ is exactly the inverse of $\mathcal{H}(\alpha)$ in $ {\rm Range}( P_{c}(H(\alpha)))$.
It suffices to prove
\begin{align}\label{111}
\frac{1}{2\pi i}{\int _\Gamma }\mathcal{H}(\alpha){\left( {\lambda  -  \mathcal{H}(\alpha)} \right)^{ - 1}}{\lambda ^{ - 1}} yd\lambda= y, \mbox{ }\forall y\in  H^2\bigcap {\rm Range}( P_{c}(H(\alpha))).
\end{align}
 Since $\mathcal{H}(\alpha)=P_{c}(H(\alpha))H(\alpha)$ only has continuous spectrum  in $X:={\rm Range}( P_{c}(H(\alpha)))$, there exists a small ball $B(0,r)$ near zero such that  $B(0,r)$ belongs to resolvent set of $\mathcal{H}(\alpha)$ in $X$. Fix $n\in \Bbb Z_+$, define the curve $\mathcal{C}_n$ to be
\begin{align*}
&\mathcal{C}_n= \mathcal{C}_\diamond\cup {\mathcal{C}_{*n}} \cup {\mathcal{C}^*_n}, \mbox{  }\mathcal{C}_{\diamond}= \partial B(0,r)\\
&\mathcal{C}_{* n} ={\mathcal{C}_{+, + }} \cup {\mathcal{C}_{+ -}} \cup {\mathcal{C}_{- -}} \cup {\mathcal{C}_{+ +}},\mbox{  }
{\mathcal{C}_{\pm\pm }} =  \{\pm a\pm t\pm i\varepsilon t: t\in[0,n]\}\\
&\mathcal{C}^*_n ={\mathcal{C}_{+ }} \cup {\mathcal{C}_{ -}},\mbox{  }{\mathcal{C}_{\pm }} =  \{t\pm i\varepsilon n: t\in[-a-n,a+n]\}.
\end{align*}
Using Cauchy integral formula on $\mathcal{C}_n$, one obtains
\begin{align*}
&\frac{1}{2\pi i} \int _\Gamma  \mathcal{H}(\alpha){\left( {\lambda  -  \mathcal{H}(\alpha)} \right)^{ - 1}}{\lambda ^{ - 1}}d\lambda  y=\lim\limits_{n\to\infty}\frac{1}{2\pi i} \int _{\mathcal{C}_{\diamond}} ...-\frac{1}{2\pi i} \int _{\mathcal{C}^{*}_{n}} ....
\end{align*}
The resolvent decay in  Lemma \ref{41} shows the integral on  ${\mathcal{C}^{*}_{n}}$ vanishes as $n\to \infty$. Recalling $B(0,r)$ belongs to resolvent set of $\mathcal{H}(\alpha)$ in $X$, the integral on  ${\mathcal{C}_{\diamond}}$ now writes
\begin{align*}
\frac{1}{2\pi i} \int_{\mathcal{C}_{\diamond}}  \mathcal{H}(\alpha){\left( {\lambda  -  \mathcal{H}(\alpha)} \right)^{ - 1}}{\lambda ^{ - 1}}d\lambda  y&= \int _{{\partial B(0,r)}} \lambda^{-1}yd\lambda+ \int_{\partial B(0,r)}{\left( {\lambda  -  \mathcal{H}(\alpha)} \right)^{ - 1}} \mathcal{H}(\alpha)yd\lambda\\
&= y.
\end{align*}
provided $y\in H^2\cap X$.  Hence, (\ref{111}) has been verified for arbitrary $y\in H^2\cap X$.

{\bf Step 2.}
Fix $0<a<a'<\tau$, let the curve $\Gamma'$ be (see Figure (c))
\begin{align*}
&\Gamma'  = {\gamma _{1, + }} \cup {\gamma _{2, + }} \cup {\gamma _{1, - }} \cup {\gamma _{2, - }},\mbox{  }{\gamma _{1, \pm }}(t) = t + a' \pm i\varepsilon t,t \ge 0,\\
&{\gamma _{2, \pm }}(t) =  - t - a' \pm i\varepsilon t,t \ge 0.
\end{align*}
In this step, we prove  for $\beta<0$
\begin{align}\label{5.1}
\int_\Gamma  {{\lambda ^\beta }} {(\lambda  - \mathcal{H}(\alpha))^{ - 1}}d\tau  = \int_{\Gamma '} {{\lambda ^\beta }} {(\lambda  - \mathcal{H}(\alpha))^{ - 1}}d\lambda.
\end{align}
Denote
{\scriptsize\begin{align*}
{\Gamma _{1,n}} &= \{  - t - a + it\varepsilon :0 \le t \le n\}  \cup \{  - t - a - it\varepsilon :0 \le t \le n\}  \\
{{\Gamma '}_{1,n}} &= \{  - t - a' + it\varepsilon :0 \le t \le n\}  \cup \{  - t - a' - it\varepsilon :0 \le t \le n\}  \\
\gamma _{1,n} &= \{  - n + it: - \varepsilon (n + a) \le t \le - \varepsilon (n + a')  \}\\
&\cup \{  - n + it: \varepsilon (n + a')\le t \le \varepsilon (n + a)\},
\end{align*}}
where $\gamma_{1,n}$ is oriented downward, $\Gamma _{1,n}$ and ${\Gamma '}_{1,n}$ are oriented anticlockwise.(see Figure (d))
Then by the analyticity of the resolvent, for $y\in X$
\begin{align*}
\int_{{\Gamma _{1,n}}} {{\lambda^\beta }} {(\lambda  - \mathcal{H}(\alpha))^{ - 1}}yd\lambda  &- \int_{{{\Gamma '}_{1,n}}} {{\lambda^\beta }} {(\lambda  - \mathcal{H}(\alpha))^{ - 1}}yd\lambda + \int_{{\gamma _{1,n}}} {{\lambda^\beta }} {(\lambda  - \mathcal{H}(\alpha))^{ - 1}}yd\lambda  = 0.
\end{align*}
Since $\Gamma$ and $\Gamma'$ are symmetric with respect to the imaginary axis, in order to prove (\ref{5.1}) it suffices to prove for $y\in X$
\begin{align}\label{100}
\mathop {\lim }\limits_{n \to \infty } \int_{{\gamma _{1,n}}} {{\lambda^\beta }} {(\lambda  - \mathcal{H}(\alpha))^{ - 1}}yd\lambda  = 0.
\end{align}
By the decay of resolvent, for $n$ sufficiently large, we have
{\small\begin{align*}
 {\left\| {\int_{{\gamma _{1,n}}} {{\lambda^\beta }} {{(\lambda  - \mathcal{H}(\alpha))}^{ - 1}} y d\lambda } \right\|_p}
&\le \int_{ - \varepsilon (n + a)}^{ - \varepsilon (n + a')} {{{\left\| {{{( - n + it - \mathcal{H}(\alpha))}^{ - 1}}y} \right\|}_p}{{\left| { - n + it} \right|}^\beta }} dt\\
&+ \int_{\varepsilon (n + a')}^{\varepsilon (n + a)} {{{\left\| {{{( - n + it - \mathcal{H}(\alpha))}^{ - 1}} y } \right\|}_p}{{\left| { - n + it} \right|}^\beta }} dt \\
&\le \int_{ - \varepsilon (n + a)}^{ - \varepsilon (n + a')} {{{\left( {{\mathop{\rm Re}\nolimits} \sqrt { \lambda \pm (it- n)} } \right)}^{ - 2}}{{\left| { - n + it} \right|}^\beta }} dt{\left\| y \right\|_p}\\
&+ \int_{\varepsilon (n + a')}^{\varepsilon (n + a)} {{{\left( {{\mathop{\rm Re}\nolimits} \sqrt {  \lambda \pm (it- n)} } \right)}^{ - 2}}{{\left| { - n + it} \right|}^\beta }} dt{\left\| y \right\|_p}.
\end{align*}}
Since as $n\to \infty$, $\arg \left( {  \lambda \pm (it- n)} \right) \sim  \pm \arctan \varepsilon$ , it is easy to see ${\left( {{\mathop{\rm Re}\nolimits} \sqrt {  \lambda \pm (it- n)} } \right)^{ - 2}} \sim c{\left| n \right|^{ - 1}}$, which immediately leads to (\ref{100}). Thus (\ref{5.1}) follows.

{\bf{Step 3}}. We prove for each $\gamma,\beta<0$,
\begin{align}\label{last}
[[\mathcal{H}(\alpha)]]^{\gamma+\beta}=[[\mathcal{H}(\alpha)]]^{\gamma}[[\mathcal{H}(\alpha)]]^{\beta}.
\end{align}
By (\ref{5.1}), for $y\in X$ one has
{\small\begin{align*}
&{[[\mathcal{H}(\alpha)]]^\gamma }{[[\mathcal{H}(\alpha)]]^\beta }y = \frac{1}{(2\pi i)^2}\int_\Gamma  {{\mu^{\gamma} }{{(\tau  - \mathcal{H}(\alpha))}^{ - 1}}} d\mu \int_{\Gamma '} {{\lambda ^\beta }{{(\lambda  - \mathcal{H}(\alpha))}^{ - 1}}y} d\lambda\\
 &=\frac{1}{(2\pi i)^2} \int_\Gamma  {\int_{\Gamma '} {{\lambda ^\beta }{\mu^{\gamma} }{{(\mu  - \mathcal{H}(\alpha))}^{ - 1}}{{(\lambda  - \mathcal{H}(\alpha))}^{ - 1}}y} d\lambda d\mu }  \\
&=\frac{1}{(2\pi i)^2} \int_\Gamma  (\mu  - \mathcal{H}(\alpha))^{ - 1}d\mu \int_{\Gamma '} \frac{\lambda ^\beta \mu^{\gamma}}{\mu-\lambda}yd\lambda-\frac{1}{(2\pi i)^2}\int_{\Gamma '} (\lambda  - \mathcal{H}(\alpha))^{ - 1} d\lambda\int_\Gamma  \frac{\lambda ^\beta \mu^{\gamma} }{\mu-\lambda}y   d\mu
\end{align*}}
which by Cauchy formula equals
\begin{align*}
&\frac{1}{2\pi i}\int_{\Gamma '} {{{(\lambda  - \mathcal{H}(\alpha))}^{ - 1}}{\lambda ^{\beta  + \gamma}} yd\lambda }={[[\mathcal{H}(\alpha)]]^{\gamma  +\beta}} y.
\end{align*}
The following two steps are standard (see Pazy \cite{AP}), but for completeness, we give a sketch.

{\bf{Step 4}}. In this step, we define $[[\mathcal{H}(\alpha)]]^{\gamma}$ for $\gamma\in {\Bbb R}$, and extend (\ref{last}) to $\gamma,\beta\in {\Bbb R}$.\\
From Step 1, we see $[[\mathcal{H}(\alpha)]]^{-1}$ is one to one in $X$. Therefore, for any integer $n\ge1$, $[[\mathcal{H}(\alpha)]]^{-n}$ is one to one. Let $\gamma<0$,
suppose $[[\mathcal{H}(\alpha)]]^{\gamma}y=0$ and take $n\ge |\alpha|$, then from Step 2, we have $[[\mathcal{H}(\alpha)]]^{-n}y= [[\mathcal{H}(\alpha)]]^{-n-\gamma}[[\mathcal{H}(\alpha)]]^{\gamma}y=0$.
 This implies $y=0$. Thus $[[\mathcal{H}(\alpha)]]^{\gamma}$ is one-to-one and we can define $[[\mathcal{H}(\alpha)]]^{-\gamma}$ as $([[\mathcal{H}(\alpha)]]^{\gamma})^{-1}$.

 We claim $[[\mathcal{H}(\alpha)]]^\gamma[[\mathcal{H}(\alpha)]]^\beta y=[[\mathcal{H}(\alpha)]]^{\gamma+\beta}y$ for $\gamma,\beta\in {\Bbb R}$ and $y\in C^{\infty}_c(\Bbb R;\Bbb C^2)\cap X$. Indeed, for example, when $\gamma>0$, $\beta<0$, $\gamma+\beta>0$, for $y\in C^{\infty}_c(\Bbb R;\Bbb C^2)\cap X$, in order to verify $[[\mathcal{H}(\alpha)]]^{\gamma+\beta}y=[[\mathcal{H}(\alpha)]]^{\gamma}[[\mathcal{H}(\alpha)]]^{\beta}y$, it suffices to prove $[[\mathcal{H}(\alpha)]]^{\beta}y=[[\mathcal{H}(\alpha)]]^{-\gamma}[[H(\alpha)]]^{\gamma+\beta}y$. Let $y=[[\mathcal{H}(\alpha)]]^{\gamma+\beta}y$, then it is equivalent to prove $[[\mathcal{H}(\alpha)]]^{-\gamma}y=[[\mathcal{H}(\alpha)]]^{\beta}[[\mathcal{H}(\alpha)]]^{-\gamma-\beta}y$, which follows immediately from Step 2.

{\bf{Step 5}}. In this step, we finish our proof.
By Step 1, $[[\mathcal{H}(\alpha)]]^1=H(\alpha)$ in $H^2\cap X$. From Step 4, $[[\mathcal{H}(\alpha)]]^{\frac{s}{2}}=[[\mathcal{H}(\alpha)]]^1[[\mathcal{H}(\alpha)]]^{\frac{s}{2}-1}
=\mathcal{H}(\alpha)[[\mathcal{H}(\alpha)]]^{\frac{s}{2}-1}$. Thus since $\frac{s}{2}-1<0$, by the definition of $\mathcal{H}(\alpha)^{\frac{s}{2}}$, we have $[[\mathcal{H}(\alpha)]]^{\frac{s}{2}}=\mathcal{H}(\alpha)^{\frac{s}{2}}$ in $H^2\cap X$. Since $[[\mathcal{H}(\alpha)]]^{\frac{s}{2}}[[\mathcal{H}(\alpha)]]^{-\frac{s}{2}}=I$, we obtain our lemma.
\end{proof}

\begin{Lemma}\label{22}
Recall $X:={\rm Range}(P_c(H(\alpha)))$.
For $2\le p\le\infty$, $f\in   X\cap H^2$, $0<s<2$, we have
\begin{align}\label{ppk}
{\left\| {{\mathcal{H}(\alpha)^{\frac{s}{2}}}f - \mathcal{H}_0^{\frac{s}{2}}f} \right\|_{L^2_x}} \le C\left\| f \right\|_{L^p_x},
\end{align}
\end{Lemma}
\begin{proof}
(\ref{20}) gives that in the space $L(H^2\cap X; L^2)$ there holds
 \begin{align*}
\mathcal{H}(\alpha)^{\frac{s}{2}} =&\int\limits_\Gamma \mathcal{H}(\alpha) {{\lambda ^{ -1+ \frac{s}{2}}}{{(\lambda  - {\mathcal{H}(\alpha)})}^{ - 1}}d\lambda }
=\int\limits_\Gamma  \mathcal{H}(\alpha){{\lambda ^{ -1+ \frac{s}{2}}}\left( {{{(\lambda  - {\mathcal{H}(\alpha)})}^{ - 1}} - {{(\lambda  -{H}_0 )}^{ - 1}}} \right)d\lambda } \\
&+ \mathcal{H}^{\frac{s}{2}}_0 + V\int\limits_\Gamma  {{\lambda ^{ -1+ \frac{s}{2}}}{{(\lambda  - {H}_0 )}^{ - 1}}d\lambda }.
\end{align*}
It is easy to see for $y\in X$ ($\mathcal{H}(\alpha)=H(\alpha)$ in $X$)
\begin{align*}
{(\lambda  -\mathcal{H}(\alpha))^{ - 1}}y - {(\lambda  -{H}_0 )^{ - 1}}y =  {(\lambda - {\mathcal{H}(\alpha)})^{ - 1}}V{(\lambda  -{H}_0 )^{ - 1}}y.
\end{align*}
Then we have
 \begin{align*}
\mathcal{H}(\alpha)^{\frac{s}{2}} =&- \int\limits_\Gamma  {(\lambda  - {\mathcal{H}(\alpha)}){\lambda ^{ -1+\frac{s}{2}}}{{(\lambda  - {\mathcal{H}(\alpha)})}^{ - 1}}V{{(\lambda  - H_0 )}^{ - 1}}d\lambda } + V\int\limits_\Gamma  {{\lambda ^{ -1+ \frac{s}{2}}}{{(\lambda  -{H}_0)}^{ - 1}}d\lambda }  \\
&+\int\limits_\Gamma {{\lambda ^{\frac{s}{2}}}{{(\lambda  - {\mathcal{H}(\alpha)})}^{ - 1}}V{{(\lambda  -{H}_0 )}^{ - 1}}d\lambda }+\mathcal{H}^{\frac{s}{2}}_0 \\
=&\int\limits_{\Gamma}  {{\lambda ^{\frac{s}{2}}}{{(\lambda  - {\mathcal{H}(\alpha)})}^{ - 1}}V{{(\lambda  -{H}_0)}^{ - 1}}d\lambda} +\mathcal{H}^{\frac{s}{2}}_0.
\end{align*}
Define
{\small \begin{align*}
&{\Gamma _1} = \Gamma  \cap \left\{ {{{({\Re}\sqrt {\tau-\lambda} )}^2} \ge 4{c_1}{{\left\| V \right\|}_\infty },\left. {{{({\Re}\sqrt { \lambda+\tau } )}^2} \ge 4{c_1}{{\left\| V \right\|}_\infty }} \right\}} \right. \\
&{\Gamma _2} = \Gamma  \cap \left\{ {{{({\Re}\sqrt {\tau-\lambda} )}^2} \le 4{c_1}{{\left\| V \right\|}_\infty },\left. {{{({\Re}\sqrt { \lambda+\tau } )}^2} \ge 4{c_1}{{\left\| V \right\|}_\infty }} \right\}} \right.\\
&{\Gamma _3} = \Gamma  \cap \left\{ {{{({\Re}\sqrt {\tau-\lambda} )}^2} \ge 4{c_1}{{\left\| V \right\|}_\infty },\left. {{{({\Re}\sqrt { \lambda+\tau } )}^2} \le 4{c_1}{{\left\| V \right\|}_\infty }} \right\}} \right. \\
&{\Gamma _4} = \Gamma  \cap \left\{ {{{({\Re}\sqrt {\tau-\lambda} )}^2} \le 4{c_1}{{\left\| V \right\|}_\infty },\left. {{{({\Re}\sqrt { \lambda+\tau } )}^2} \le 4{c_1}{{\left\| V \right\|}_\infty }} \right\}} \right..
\end{align*}}
Then from Lemma \ref{41} when $\frac{1}{q}+\frac{1}{p}=\frac{1}{2}$, one deduces
{\small \begin{align*}
&{\left\| {\int\limits_{{\Gamma _1}} {{\lambda ^{\frac{s}{2}}}{{(\lambda  - H(\alpha))}^{ - 1}}V{{(\lambda  - {H_0})}^{ - 1}}fd\lambda } } \right\|_2} \\
&\le \int\limits_{{\Gamma _1}} {{{\left| \lambda  \right|}^{\frac{s}{2}}}{{\left\| {{{(\lambda  - H(\alpha))}^{ - 1}}} \right\|}_{2 \to 2}}{{\left\| V \right\|}_q}{{\left\| {{{(\lambda  - {H_0})}^{ - 1}}} \right\|}_{p \to p}}{{\left\| f \right\|}_p}d\lambda }
\end{align*}}
which by Lemma \ref{45} is further dominated by
\begin{align}
\int\limits_{\Gamma_1}  {{{\left| \lambda  \right|}^{\frac{s}{2}}}[{{({\Re}\sqrt {\tau-\lambda} )}^{ - 2}} + {{({\Re}\sqrt { \lambda+\tau } )}^{ - 2}}]^2{{\left\| V \right\|}_q}{{\left\| f \right\|}_p}d\lambda }.\label{42}
\end{align}
The only singular point for (\ref{42}) is infinity. But at infinity, the integrand behaves like
${{{\left| \lambda  \right|}^{\frac{s}{2} - 2}}}$, which is integrable.
Notice that
${({\Re}\sqrt {\tau-\lambda} )^{ - 2}}\sim c|\lambda|^{-1}$ , ${({\Re}\sqrt { \lambda+\tau } )^{ - 2}} \sim c \left| \lambda  \right|^{-1}$,
as $\lambda\to\infty$ in $\Gamma$, the curve $\Gamma_2$, $\Gamma_3$, $\Gamma_4$ are actually bounded. Hence the remaining three terms can be estimated by
\begin{align*}
&\int\limits_{\left\{ {\lambda :\left| \lambda  \right| \le R} \right\} \cap \Gamma } {{{\left| \lambda  \right|}^{\frac{s}{2}}}{{\left\| {{{\left( {\lambda  - H(\alpha)} \right)}^{ - 1}}} \right\|}_{2 \to 2}}{{\left\| V \right\|}_q}{{\left\| {{{\left( {\lambda  - {H_0}} \right)}^{ - 1}}} \right\|}_{p \to p}}{{\left\| f \right\|}_p}d\lambda }  \le C{\left\| V \right\|_q}{\left\| f \right\|_p}.
\end{align*}
\end{proof}

The same arguments give
\begin{Lemma}\label{ty678}
For any $0<s<2$ and $\forall f\in H^2$, one has
\begin{align}\label{109b1}
\|{\mathcal{H}(\alpha_1)^{\frac{s}{2}}}f-{\mathcal{H}(\alpha)^{\frac{s}{2}}}f\|_{2}  \le C(|\alpha-\alpha_1|)|\left\| f \right\|_2,
\end{align}
\end{Lemma}
\begin{proof}
If $f\in {\rm Range} \left(I-P_c(H(\alpha))\right)$ or   $f\in{\rm Range}\left(I-P_c(H(\alpha_1))\right)$, (\ref{109b1}) holds naturally because the generalized eigenfunctions $\{\xi_i\}$ are Schwartz functions and depends smoothly on parameter $\alpha$. If $f\in {\rm Range} \left(P_c(H(\alpha))\right)\cap {\rm Range} \left(P_c(H(\alpha_1))\right)$, using same arguments as Lemma \ref{22} and the smooth dependence of potential $V$ on $\alpha$ gives (\ref{109b1}).
\end{proof}

\begin{Lemma}\label{46}
For $1/2<s<2$, $2\le p\le\infty$, $\forall f\in H^2$, it holds
$$\|f\|_p+\|f\|_{H^s}\le C\|\mathcal{H}_0^{\frac{s}{2}}f\|_2$$
\end{Lemma}
\begin{proof}
By Lemma \ref{41}, the integrand in (\ref{109b}) is absolutely integrable in $H^2_x$. Hence, the Fourier transform denoted by $\mathcal{F}$ on $\Bbb R$ can go across the integral symbol and act directly on the integrand in (\ref{109b}) by Plancherel theorem.  Thus Cauchy integral formula yields
{\footnotesize\begin{align*}
& \mathcal{F}\left[ {\int_\Gamma  {{\lambda ^{\frac{s}{2} - 1}}} {{(\lambda  - {H_0})}^{ - 1}}{{H}_0}fd\lambda } \right](k) \\
 & = \left( \begin{array}{l}
 \tau  + {\left| k \right|^2} \\
  \\
 \end{array} \right.\left. \begin{array}{l}
  \\
  - \tau  - {\left| k \right|^2} \\
 \end{array} \right)
 \int_\Gamma  {{\lambda ^{\frac{s}{2} - 1}}} \left( \begin{array}{l}
 {\left[ {\lambda  - (\tau  + {{\left| k \right|}^2})} \right]^{ - 1}} \\
  \\
 \end{array} \right.\left. \begin{array}{l}
  \\
 {\left[ {\lambda  + \tau  + {{\left| k \right|}^2}} \right]^{ - 1}} \\
 \end{array} \right)\mathcal{F}(f)(k)d\lambda  \\
 &= \left( \begin{array}{l}
 {\left( {\tau  + {{\left| k \right|}^2}} \right)^{\frac{s}{2}}} \\
  \\
 \end{array} \right.\left. \begin{array}{l}
  \\
 {\left( { - \tau  - {{\left| k \right|}^2}} \right)^{\frac{s}{2}}} \\
 \end{array} \right)\mathcal{F}(f)(k). \\
 \end{align*}}
Thus $H_0^{\frac{s}{2}}$ is roughly $(\tau-\Delta)^{\frac{s}{2}}$. Then standard Sobolev embedding yields ${\left\| f \right\|_p} \le C{\left\| {\mathcal{H}_0^{\frac{s}{2}}f} \right\|_2}$.
\end{proof}

\begin{Lemma}\label{31a}
For $f\in H^2$ in {\bf the continuous spectrum space} of $H(\alpha)$, $1/2<s<2$, $2\le p\le \infty$, we have
$$\|f\|_p\le C\|\mathcal{H}(\alpha)^{\frac{s}{2}}f\|_2.$$
\end{Lemma}
\begin{proof}
It suffices to prove
${\left\| {{{[[\mathcal{H}(\alpha)]]}^{ - \frac{s}{2}}}f} \right\|_p} \le C{\left\| f \right\|_2},$ due to Lemma \ref{51a}.
Similar arguments as Lemma \ref{22} yield,
\begin{align*}
&{\left\| {{{[[H(\alpha)]]}^{ - \frac{s}{2}}}f - {{[[{\mathcal{H}_0}]]}^{ - \frac{s}{2}}}f} \right\|_p} \le C{\left\| {\int_\Gamma  {{\lambda ^{ - \frac{s}{2} + 1}}{{(\lambda  - {{H}_0})}^{ - 1}}V{{(\lambda  - \mathcal{H}(\alpha))}^{ - 1}}f} d\lambda } \right\|_p} \\
&\le C{\int_{\Gamma_1} {{{\left| \lambda  \right|}^{ - \frac{s}{2} + 1}}{{\left\| {{{(\lambda  - {H_0})}^{ - 1}}} \right\|}_{2 \to p}}{{\left\| V \right\|}_\infty }\left\| {{{(\lambda  - \mathcal{H}(\alpha))}^{ - 1}}} \right\|} _{2 \to 2}}{\left\| f \right\|_2}d\lambda+II,
\end{align*}
where $II$ denotes the remaining $\Gamma_2$, $\Gamma_3$, $\Gamma_4$ parts.
Applying Lemma \ref{41}, Lemma \ref{45} and Lemma \ref{46}, using similar arguments as Lemma \ref{22}, we obtain Lemma \ref{31a}.
\end{proof}

\section{Orthogonal conditions}
For vector valued function $\phi$,
define
\begin{align*}
&{\left| {{J_V(\alpha)}} \right|^s}\phi \\
&= \left( {\begin{array}{*{16}{c}}
   {M(t+h)}  \\
   {}  \\
\end{array}} \right.\left. {\begin{array}{*{16}{c}}
   {}  \\
   {M( - t-h)}  \\
\end{array}} \right){\left( {{(t+h)^2}\mathcal{H}(\alpha)} \right)^{\frac{s}{2}}}\left( {\begin{array}{*{16}{c}}
   {M( - t-h)}  \\
   {}  \\
\end{array}} \right.\left. {\begin{array}{*{16}{c}}
   {}  \\
   {M(t+h)}  \\
\end{array}} \right)\phi,
\end{align*}
where $M(t)=e^{i|x|^2/4t}$, $h>0$.
$\mathcal{H}(\alpha)^{\frac{s}{2}}$ has been defined in Section 2. Here, we recall the fact that $\mathcal{H}(\alpha)^{\frac{s}{2}}\xi_i=0$ and $(\mathcal{H}^*(\alpha))^{\frac{s}{2}}\theta_3=\theta_3\mathcal{H}(\alpha)^{\frac{s}{2}}$.

Assume that the solution to (1.1) is of the following form:
\begin{align}\label{p12}
u(x,t) = w(x;\sigma (t)) + \chi (x,t),w(x; \sigma (t)) = e^{- i\beta (t)}\varphi (x;\alpha (t)),
\end{align}
where $\sigma (t)=(\beta(t),\omega(t),0,0)$ is not a solution of (\ref{70}) in general.
Define $\chi (x,t)=e^{-i\beta(t)}f(x,t)$, $\mathbf{f}=(f,\bar{f})^t$ and
$$U(t)=\left( \begin{array}{l}
 M(t) \\
  \\
 \end{array} \right.\left. \begin{array}{l}
  \\
 M(-t) \\
 \end{array} \right).
$$

\begin{Lemma}
For any $\phi\in H^2$, one has
\begin{align}\label{key}
\left\langle {{{|J_V(\alpha)|^s}}{\phi},U(t+h){\theta _3}{\xi _i(\alpha)}} \right\rangle  = 0,
\end{align}
for $i=1,2$.
\end{Lemma}
\begin{proof}
Since $(\mathcal{H}^*(\alpha))^{\frac{s}{2}}\theta_3=\theta_3\mathcal{H}(\alpha)^{\frac{s}{2}}$, it holds
\begin{align*}
 \left\langle { {{|J_V(\alpha)|^s}}{\phi},U(t+h){\theta _3}{\xi _i}} \right\rangle
  &= (t+h)^s\left\langle {U(-t-h){\phi},{{\left( {{\mathcal{H}^*}({\alpha})} \right)}^{\frac{s}{2}}}{\theta _3}{\xi _i}} \right\rangle  \\
  &= (t+h)^s\left\langle {U(-t-h){\phi},{\theta _3}{{\left( {\mathcal{H}({\alpha})} \right)}^{\frac{s}{2}}}{\xi _i}} \right\rangle.
\end{align*}
Since $\mathcal{H}(\alpha)^{\frac{s}{2}}\xi_i=0$(by Lemma \ref{xinhui}), we obtain
\begin{align*}
&\left\langle {U(-t-h){\phi},{\theta _3}{{\left( {\mathcal{H}({\alpha})} \right)}^{\frac{s}{2}}}{\xi _i}} \right\rangle  = 0.
\end{align*}
\end{proof}

We impose two orthogonal conditions to $f$: for $i=1,2$
\begin{align}\label{777}
\left\langle {U(-t-h)\mathbf{f},{\theta _3}{\xi _i(\alpha(t))}} \right\rangle  = 0.
\end{align}
The existence of $\sigma(t)$ follows by the following lemma.

\begin{Lemma}\label{ghdftry}
If $h$ is sufficiently large, $\|\chi(0,x)\|_{L^2}$ is sufficiently small, then there exists $\sigma(t)=(\beta(t),\omega(t),0,0)$ such that (\ref{777}) holds.
\end{Lemma}
\begin{proof}
We first prove it for $t=0$, namely there exist appropriate $\alpha$ and $\beta$ such that
$$\left\langle {\left( \begin{array}{l}
 u_0(x) - {e^{ - i\beta }}\varphi (x;\alpha ) \\
 \overline {u_0(x)}  - {e^{i\beta }}\varphi (x;\alpha ) \\
 \end{array} \right),\left( \begin{array}{l}
 {e^{ - i\beta }} \\
  \\
 \end{array} \right.\left. \begin{array}{l}
  \\
 {e^{i\beta }} \\
 \end{array} \right)U(-h){\theta _3}{\xi _i}(\alpha )} \right\rangle  = 0.
 $$
This solvability is a consequence of the non-singularity of the corresponding Jacobian. Indeed, since $|M(t)-1|\le C\frac{|x|^2}{t}$, $\chi(0,x)$ is small in $L^2$, for $h$ sufficiently large, the leading term of the Jacobian is
$$\left( \begin{array}{l}
 e \\
 0 \\
 \end{array} \right.\left. \begin{array}{l}
 0 \\
 s \\
 \end{array} \right),
$$
where $e =  - i\frac{d}{{d\alpha }}\left\| {\varphi (x;\alpha )} \right\|_2^2$, $s =  - \frac{2}{\alpha }i\frac{d}{{d\alpha }}\left\| {\varphi (x;\alpha )} \right\|_2^2.$
At $\alpha_0$, $e$ and $s$ are nonzero by the assumption in Theorem 1.1. Thus we have proved our lemma when $t=0$.\\
Second we show the existence of $\sigma(t)$ for $t>0$, but this follows by the standard arguments and the orbital stability, see \cite{BP} for instance.
\end{proof}

Lemma \ref{31a} and the definition of $|J_V(\alpha)|^su$ immediately yield
\begin{Corollary}\label{54}
For $\phi\in H^2$, $\frac{1}{2}<s<2$, $2\le p\le\infty$, we have
$$\|P_{c}(H(\alpha))(U(-t-h)\phi)\|_{p}\le t^{-s}\||J_V(\alpha)|^s\phi\|_{L^2_x}.
$$
provided that the RHS is finite.
\end{Corollary}
\begin{proof}
Decompose $U(-t-h)\phi$ into
\begin{align*}
U(-t-h)\phi=P_{c}(H(\alpha))[U(-t-h)\phi]+\mu_1\xi_1(\alpha)+\mu_2\xi_2(\alpha).
\end{align*}
Applying Lemma \ref{xinhui} shows the discrete part $\sum \mu_i\xi_i$ vanishes in $|J_V(\alpha)|^s\phi$. And then using Lemma \ref{31a} yields our result.
\end{proof}

\section{Set up of the bootstrap argument}

Let ${\bf f}=(f,\bar{f})^t$ be the remainder term in the decomposition (\ref{p12}) such that (\ref{777}) holds. (see Lemma  \ref{ghdftry}) For any given $t_1>0$, define
\begin{align}\label{phk89n}
{\mathcal{M}}_{t_1} = \mathop {\sup }\limits_{0\le\tau  \le t_1} ({\left\| f \right\|_m} + {\left\| f \right\|_n} + {\left\| f \right\|_\infty })
{\left\langle \tau  \right\rangle ^s}.
\end{align}
Let $\mathcal{A}\subset [0,\infty)$ be the set of $t_1$ such that $\mathcal{M}_{t_1}$ is sufficiently small, i.e.,
\begin{align}
\mathcal{A}=\{t_1\in[0,\infty):{\mathcal{M}_{t_1}}<\kappa\ll1\},
\end{align}
where $\kappa>0$ is sufficiently small to be determined in Section 6. $\mathcal{A}$ is nonempty due to (\ref{p0kinbgh}) and Sobolev embedding. Moreover, $\mathcal{A}$ is closed by the global well-posedness theory and Sobolev embeddings.

\subsection{Modulation equation}

In this subsection, we study the modulation equation.
Define $\beta(t)=\int^t_0 \omega(\lambda)d\lambda +\gamma$, then we rewrite (1.1) in terms of $f$:
$$if_t=L(\alpha(t))f+N(\varphi,f)-\gamma'f -\gamma'\varphi+i\omega'\frac{2}{\alpha}\varphi_\alpha,$$
where $N(\varphi,f)=F(|\varphi+f|^2)(\varphi+f)-F(\varphi^2)\varphi-F(\varphi^2)f-F'(\varphi^2)\varphi^2(f+\bar{f}).$
By (\ref{777}), we have
\begin{align}\label{6}
{\rm{Im}}\left\langle {f,M(t+h){v_i}} \right\rangle  = 0.
\end{align}

\begin{Lemma}\label{a1}
If $t_1\in\mathcal{A}$, then for $t\in[0,t_1]$, $\|f(t)\|_{\infty}$ is sufficiently small, and
$$|\gamma '(t)| + |\omega '(t)| \le C(\left\| f \right\|_{_\infty }^2 + \left\| f \right\|_\infty ^{m - 1}{\left\| f \right\|_m} + \left\| f \right\|_\infty ^{n - 1}{\left\| f \right\|_n} + \frac{1}{{{t^2}}}{\left\| f \right\|_\infty } + \frac{1}{t}{\left\| f \right\|_\infty }).
$$
\end{Lemma}
\begin{proof}
Since $t_1\in\mathcal{A}$ and $t\in[0,t_1]$, $\|f(t)\|_{\infty}$ is sufficiently small by (\ref{phk89n}).
Differentiating (\ref{6}) with respect to $t$, we obtain
\begin{align}\label{7}
&{\rm{Im}}\left\langle {i\gamma '\varphi  + \omega '\frac{2}{\alpha }{\varphi _\alpha },M(t+h){v_i}} \right\rangle  + {\rm{Im}}\left\langle {f,M(t+h){\partial _\alpha }{v_i}} \right\rangle \omega '\frac{2}{\alpha } \nonumber\\
&+ {\rm{Im}}\left\langle { - iL(\alpha (t))f,M(t+h){v_i}} \right\rangle+ {\rm{Im}}\left\langle { - iN(\varphi ,f),M(t+h){v_i}} \right\rangle\nonumber\\
&+ {\rm{Im}}\left\langle {i\gamma 'f,M(t+h){v_i}} \right\rangle  + {\rm{Im}}\left\langle {f,\frac{d}{{dt}}M(t+h){v_i}} \right\rangle  = 0.
\end{align}
Direct calculations give
\begin{align}\label{8}
\left| {\left\langle {f,\frac{d}{{dt}}M(t+h){v_i}} \right\rangle } \right| \le {t^{ - 2}}{\left\| f \right\|_{\infty}}.
\end{align}
By the orthogonal condition (\ref{777}), the identity $H(\alpha)\xi_2=i\xi_1$ and the obvious commutator inequality, we have
\begin{align}
&\left| {{\rm{Im}}\left\langle { - iL(\alpha (t))f,M(t+h){v_i}} \right\rangle } \right|=\left| {\left\langle {H(\alpha (t)){\bf f},U(t+h){\theta _3}{\xi_i}} \right\rangle } \right|\nonumber\\
&= \left| {\left\langle {{\bf f},U(t+h){H}^*(\alpha (t)){\theta _3}{\xi_i}} \right\rangle } \right| + \left| {\left\langle {{\bf f},[U(t+h),H^*(\alpha (t))]{\theta _3}{\xi_i}} \right\rangle } \right|\nonumber\\
&= \left| {\left\langle {{\bf f},U(t+h){\theta _3}H(\alpha (t)){\xi_i}} \right\rangle } \right| + \left| {\left\langle {{\bf f},[U(t+h),H^*(\alpha (t))]{\theta _3}{\xi_i}} \right\rangle } \right|\nonumber\\
&= \left| {\left\langle {{\bf f},[U(t+h),H^*(\alpha (t))]{\theta _3}{\xi_i}} \right\rangle } \right|
\le {(t+h)^{ - 1}}{\left\| f \right\|_{\infty}}.\label{poiuj}
\end{align}
Combining (\ref{poiuj}) with (\ref{7}) and (\ref{8}) yields
\begin{align}\label{9}
\left\| {\left( {\begin{array}{*{20}{c}}
   0  \\
   { - e}  \\
\end{array}} \right.\left. {\begin{array}{*{20}{c}}
   e \\
   0  \\
\end{array}} \right)\left( {\begin{array}{*{20}{c}}
   {\gamma '}  \\
   {\omega '}  \\
\end{array}} \right)} \right\| \lesssim {\left\| f \right\|_{\infty}}\left\| {\left( {\begin{array}{*{20}{c}}
   {\gamma '}  \\
   {\omega '}  \\
\end{array}} \right)} \right\| + {O_2}(f) + {t^{ - 1}}{\left\| f \right\|_{\infty}} + {t^{ - 2}}{\left\| f \right\|_{\infty}},
\end{align}
where $O_2\le C(|f|^m+|\varphi|^{m-2}|f|^2+|f|^n+|\varphi|^{n-2}|f|^2)$.
Thus our lemma follows by the smallness of $\|f\|_{\infty}$.
\end{proof}

For $s=\frac{7}{4}^+$,  Lemma \ref{a1} implies
\begin{Corollary}\label{po987uj}
For any $s=\frac{7}{4}^+$, any $t_1\in\mathcal{A}$ and arbitrary $0\le t\le t_1$, we have
\begin{align}\label{zx}
\left| {\gamma '} (t)\right| + \left| {\omega '} (t)\right| \le CW(\mathcal{M}_{t_1}) {t^{ - 1 - s}}\mathcal{M}^2_{t_1},
\end{align}
where $W(x)$ is some bounded function around zero.
\end{Corollary}

\subsection{Reduction to time-independent linearized operator}

Fix any time $t_1\in \mathcal{A}$. Suppose that $\omega(t_1)=\omega_1$, $\beta(t_1)=\beta_1$. Define $\omega_1=-\frac{1}{4}\alpha^2_1$, $\gamma_1=\beta_1-\omega_1t_1$, $\Phi_1=-\omega_1t-\gamma_1$. Let $\chi=e^{i\Phi_1}g(x,t)$, then $g$ satisfies
\begin{align}\label{kiou8}
ig_t=L(\alpha_1)g+D.
\end{align}
The function $D$ is given by
\begin{align*}
 D &= {D_0} + {D_1} + {D_2} + {D_3} + {D_4}, \\
 {D_0} &= {e^{ - i\Omega }}[ - \gamma '\varphi (\alpha ) + \frac{{2i}}{\alpha }\omega '{\varphi _\alpha }], \mbox{  }\mbox{  }\Omega  = {\Phi _1} - \Phi , \\
 {D_1} &= [F({\varphi ^2}(\alpha )) + F'({\varphi ^2}(\alpha )){\varphi ^2}(\alpha ) - F({\varphi ^2}({\alpha _1})) - F'({\varphi ^2}({\alpha _1})){\varphi ^2}({\alpha _1})]g, \\
 {D_2} &= F'({\varphi ^2}(\alpha )){\varphi ^2}(\alpha )[{e^{ - 2i\Omega }} - 1]\bar g, \\
 {D_3} &= [F'({\varphi ^2}(\alpha )){\varphi ^2}(\alpha ) - F'({\varphi ^2}({\alpha _1})){\varphi ^2}({\alpha _1})]\bar g, \\
 {D_4} &= {e^{ - i\Omega }}N(\varphi (\alpha ),{e^{i\Omega }}g).
\end{align*}
The equation for $\mathbf{g}=(g,\bar{g})^t$ is
\begin{align}\label{b1}
i\mathbf{g}_t=H(\alpha_1)\mathbf{g}+\mathbf{D}.
\end{align}

For the given time $t_1\in\mathcal{A}$ denote $\alpha(t_1)$ by $\alpha_1$.
Recall
$${H(\alpha_1)} = \left( \begin{array}{l}
  - \Delta  + \frac{{\alpha _1^2}}{4} \\
  \\
 \end{array} \right.\left. \begin{array}{l}
  \\
 \Delta  - \frac{{\alpha _1^2}}{4} \\
 \end{array} \right) + \left( \begin{array}{l}
 {V_1}({\alpha _1}) \\
 {V_2}({\alpha _1}) \\
 \end{array} \right.\left. \begin{array}{l}
 {V_2}({\alpha _1}) \\
  - {V_1}({\alpha _1}) \\
 \end{array} \right).
 $$

Recall that $\bf{g}$ is the remainder in the above decomposition of $u(x,t)$ (see (\ref{kiou8})).
Moreover, let $s=\frac{7}{4}^+$, then by Corollary \ref{po987uj} and Newton-Leibnitz formula, for $t_1\in\mathcal{A}$, $t\in[0,t_1]$,
\begin{align}
\mathop {\sup }\limits_{0\le\tau  \le t_1} \left| {\alpha (t) - {\alpha _0}} \right|  &\le W(\mathcal{M}_{t_1})\mathcal{M}^2_{t_1},\nonumber\\
\left| {\alpha '} \right| &\le W({\mathcal{M}_{t_1}}){t^{ - 1 - s}}{{\mathcal{M}_{t_1}}^2}, \label{60} \\
\left| {\alpha (t) - \alpha ({t_1})} \right|+|\Omega| &\le W(\mathcal{M}_{t_1}){\mathcal{M}^2_{t_1}}\left \langle t\right \rangle^{1 - s} . \label{s1}
\end{align}
Section 1 implies the discrete spectral part of $H(\alpha)$ is spanned by $\{\xi_1(\alpha),\xi_2(\alpha)\}$ in the radial case. Moreover it is known that $\xi_i(\alpha)$ depends continuously with respect to $\alpha$.
Denote the projection to the continuous part of $H(\alpha_1)$ as $P_2$, i.e. $P_2=P_c(H(\alpha_1))$.

\begin{Lemma}\label{31}
Let $p\in[2,\infty]$, $h\gg1$, $t_1\in\mathcal{A}$, then we have for $0\le t\le t_1$
\begin{align}
{\left\| U(-t-h)\mathbf{g}(t) \right\|_p} &\lesssim {\left\| {{P_2}U(-t-h)\mathbf{g}} \right\|_p}.\label{wushu1}
\end{align}
\end{Lemma}
\begin{proof}
First of all by (\ref{s1}) and $t_1\in\mathcal{A}$, we have
\begin{align}\label{p098nhbg}
|\alpha_1-\alpha(t)|=|\alpha(t_1)-\alpha(t)|\ll1.
\end{align}
Recall $g=e^{-i\Omega}f$.
Then by (\ref{777}),
\begin{align}
0& = \left\langle {\mathbf{f},U( - t - h){\xi _i}(\alpha (t))} \right\rangle  = \left\langle {\mathbf{g},{e^{ - i\Omega }}U( - t - h){\xi _i}(\alpha (t))} \right\rangle\nonumber\\
&= \left\langle { U(-t-h)\mathbf{g},{e^{ - i\Omega }}U( - 2t -2 h){\xi _i}(\alpha_1)} \right\rangle\nonumber\\
&+\left\langle { U(-t-h)\mathbf{g},{e^{ - i\Omega }}U( - 2t - 2h)({\xi _i}(\alpha(t))-{\xi _i}(\alpha_1))} \right\rangle.\label{wushu3}
\end{align}
By (\ref{p098nhbg}) and $\xi_i$ is of exponential decay as $x\to\infty$,
\begin{align}\label{wushuyu3}
|\left\langle { U(-t-h)\mathbf{g},{e^{ - i\Omega }}U( - t - h)({\xi _i}(\alpha(t))-{\xi _i}(\alpha_1))} \right\rangle|\ll \|U(-t-h)\mathbf{g}\|_{p}.
\end{align}
Splitting $U(-t-h)\bf g$ into continuous spectral part and discrete part,
$$ U(-t-h)\mathbf{g}= {P_c}(H(\alpha_1)) U(-t-h)\mathbf{g}  + {\mu _1}{\xi _1}(\alpha_1) + {\mu _2}{\xi _2}(\alpha_1),$$
we have by (\ref{wushu3}) and (\ref{wushuyu3}) that
\begin{align*}
&{\mu _1}\left\langle {{\xi _1}(\alpha_1),{e^{ - i\Omega }}U( - t - h){\xi _i}(\alpha_1)} \right\rangle + {\mu _2}\left\langle {{\xi _2}(\alpha_1),{e^{ - i\Omega }}U( - t - h){\xi _i}(\alpha_1)} \right\rangle \\
& =  O(\nu(|\mu_1|+|\mu_2|))+O(\|{P_2} U(-t-h)\mathbf{g}\|_p).
\end{align*}
with some constant $0<\nu\ll1$.
And (\ref{s1}) gives
$$\left\| {{e^{ - i\Omega }}U( - t - h) - I} \right\| \le W({{\mathcal M}_{t_1}}){{{\mathcal M}_{t_1}}^2}{\left\langle t\right\rangle ^{1 - s}} + \frac{1}{{t + h}}{\left| x \right|^2}.
$$
Hence by the exponential decay of $\xi_i$, if ${\mathcal M}_{t_1}$ is sufficiently small and $h\gg 1$, we obtain
\begin{align*}
\|(\mu_1,\mu_2)\|\lesssim  \|P_2U( - t - h)\mathbf{g}\|_p,
\end{align*}
thus proving (\ref{wushu1}).
\end{proof}

\begin{Lemma}\label{a2}
Let $s=\frac{7}{4}^+$. If $t_1\in\mathcal{A}$, $0\le t\le t_1$, then there exists some universal constant $C>0$ such that
\begin{align}
\||{J_V(\alpha_1)}{|^s}\mathbf{g}(t)\|_2&\le C\|P_2|{J_V(\alpha_1)}|^s\mathbf{g}\|_2\nonumber\\
{\left\| \mathbf{g} \right\|_p} &\le C {(t + h)^{ - s}}{\left\| {{{\left| {{J_V}(\alpha_1)} \right|}^s}\mathbf{g}} \right\|_2}.\label{wushu2}
\end{align}
\end{Lemma}
\begin{proof}
Suppose $|J_V(\alpha_1)|^s\mathbf{g}=P_2|J_V(\alpha_1)|^s\mathbf{g}+k_1\xi_1(\alpha_1)+k_2\xi_2(\alpha_1)$, then (\ref{key}) gives
\begin{align}
0&= \left\langle {{|J_V(\alpha(t))|^s}\mathbf{g},U(t+h){\theta _3}{\xi _i}(\alpha(t) )} \right\rangle  \nonumber\\
 &= \left\langle {|{J_V(\alpha_1)}{|^s}\mathbf{g},U(t+h){\theta _3}[{\xi _i}(\alpha (t)) - {\xi _i}({\alpha _1})]} \right\rangle\label{p910i}\\
 &+\left\langle {|{J_V(\alpha_1)}{|^s}\mathbf{g},U(t+h){\theta _3}{\xi _i}({\alpha _1})} \right\rangle\nonumber\\
 &+\left\langle {(|{J_V(\alpha(t))}{|^s}-|{J_V(\alpha_1)}{|^s})\mathbf{g},U(t+h){\theta _3}{\xi _i}({\alpha(t)})} \right\rangle\label{p9i}
\end{align}
{\bf Step 1.}
(\ref{p9i}) is bounded by $\|(|{J_V(\alpha(t))}{|^s}-|{J_V(\alpha_1)}{|^s})\mathbf{g}\|_2$, which by Lemma \ref{ty678} has an upper bound of $|\alpha_1-\alpha|\|\mathbf{g}\|_2$. Hence Lemma \ref{31} shows (\ref{p9i}) is bounded by $\nu{\left\| {{P_2}U(-t-h)\mathbf{g}} \right\|_2}$ with some constant $0<\nu\ll1$. (notice the matrix $U(t)$ is unitary ) Then Corollary \ref{54} implies
\begin{align*}
|(\ref{p9i})| \le  \nu\||{J_V(\alpha_1)}|^s\mathbf{g}\|_2.
\end{align*}
{\bf Step 2.}
Meanwhile, (\ref{p910i}) is bounded by
\begin{align*}
\left\langle {{{\left| {{J_V(\alpha_1)}} \right|}^s}\mathbf{g},{\theta _3}{\xi _i}({\alpha _1})} \right\rangle  + O\left( {\left| {{\alpha _1} - \alpha } \right|} \right){\left\| {{{\left| {{J_V(\alpha_1)}} \right|}^s}\mathbf{g}} \right\|_2} + \|( U(t+h)-I){\xi _i}\|_{L^{\infty}_x}\left\|\left| J_V(\alpha_1) \right|^s\mathbf{g} \right\|_2.
\end{align*}
Therefore substituting the spectrum decomposition of $|J_V|^s\mathbf{g}$ into the above formula yields
\begin{align*}
\left\| {({k_1},{k_2})} \right\| &\lesssim {\left\| {{P_2}{{\left| {{J_V(\alpha_1)}} \right|}^s}\mathbf{g}} \right\|_2} + O\left( {\left| {{\alpha _1} - \alpha } \right|} \right){\left\| {{{\left| {{J_V(\alpha_1)}} \right|}^s}\mathbf{g}} \right\|_2}+ h^{-1}{\left\| {{{\left| {{J_V}}(\alpha_1) \right|}^s}\mathbf{g}} \right\|_2} \\
&\lesssim {\left\| {{P_2}{{\left| {{J_V}(\alpha_1)} \right|}^s}\mathbf{g}} \right\|_2} + O\left( {\left| {{\alpha _1} - \alpha } \right|} \right)\left\| {({k_1},{k_2})} \right\|
+ h^{-1}{\left\| {{P_2{\left| {{J_V}(\alpha_1)} \right|}^s}\mathbf{g}} \right\|_2}\\
&+h^{-1}{\left\| (k_1,k_2)\right\|}.
\end{align*}
Thus by (\ref{p098nhbg}) and $h\gg 1$, we have
$$\left\|(k_1,k_2) \right\| \le C{\left\| {{P_2}{{\left| {{J_V}}(\alpha_1) \right|}^s}\mathbf{g}} \right\|_2}.$$
Finally, notice that (\ref{wushu2}) is a direct corollary of (\ref{wushu1}) and Corollary \ref{54}.
\end{proof}

Therefore, Corollary \ref{54} and Lemma \ref{a2} yield
\begin{Corollary}
Suppose that $t_1\in\mathcal{A}$, $s=\frac{7}{4}^+$, then
\begin{align}\label{kill}
\mathcal{M}_{t_1}\le C\|P_2|J_V(\alpha_1)|^s\mathbf{g}\|_{L_t^\infty[0,t_1] L_x^2}.
\end{align}
\end{Corollary}

\section{The equation for $|J_V(\alpha_1)|^sg$}
We will prove Theorem 1.1 by bootstrap. As a preparation, we derive the equation of $|J_V(\alpha_1)|^sg$ first.

\label{kill}
Denote $${\partial _t}{E_2} = \left( \begin{array}{l}
 {\partial _t} \\
  \\
 \end{array} \right.\left. \begin{array}{l}
  \\
 {\partial _t} \\
 \end{array} \right),
 $$
then
direct calculations imply
\begin{Lemma}\label{50}
{\small $$\begin{array}{l}
 \left[ {i{\partial _t}{E_2} - \left( \begin{array}{l}
  - \Delta  + \tau \\
  \\
 \end{array} \right.\left. \begin{array}{l}
  \\
 \Delta  -\tau \\
 \end{array} \right),\left( \begin{array}{l}
 M(t+h) \\
  \\
 \end{array} \right.\left. \begin{array}{l}
  \\
 M( -t-h) \
 \end{array} \right)} \right] \\
  = \left( \begin{array}{l}
 M(t+h)(\frac{i}{{2(t+h)}} + \frac{{ix \cdot \nabla }}{t+h}) \\
  \\
 \end{array} \right.\left. \begin{array}{l}
  \\
 M(-t-h)(\frac{i}{{2(t+h)}} + \frac{{ix \cdot \nabla }}{t+h}) \\
 \end{array} \right). \\
 \end{array}
 $$}
\end{Lemma}

\begin{Lemma}\label{huziji9}
{\small $$
 \left[ {i{\partial _t}{E_2} - \left( \begin{array}{l}
  - \Delta  + \tau\\
  \\
 \end{array} \right.\left. \begin{array}{l}
  \\
 \Delta  - \tau \\
 \end{array} \right),\left( \begin{array}{l}
 M( - t-h) \\
  \\
 \end{array} \right.\left. \begin{array}{l}
  \\
 M(t+h) \\
 \end{array} \right)} \right] \\
  =  \mathfrak{M}
 $$}
\end{Lemma}
where we denote
 {\small \begin{align*}
 \mathfrak{M}:=
 \left( \begin{array}{l}
 M(-t-h)( - \frac{i}{{2(t+h)}} - \frac{{ix \cdot \nabla }}{t+h} - \frac{{{x^2}}}{{2{(t+h)^2}}}) \\
  \\
 \end{array} \right.\left. \begin{array}{l}
  \\
 M(t+h)( - \frac{i}{{2(t+h)}} - \frac{{ix \cdot \nabla }}{t+h} + \frac{{{x^2}}}{{2{(t+h)^2}}}) \\
 \end{array} \right).
  \end{align*}}

\noindent For the potential term, we have
\begin{Lemma}\label{51}
$$\left[ {V,U(t)} \right] = \left( {\begin{array}{*{20}{c}}
   0  \\
   {{V_2}[M( - t-h) - M(t+h)]}  \\
\end{array}} \right.\left. {\begin{array}{*{20}{c}}
   {-{V_2}[M( - t-h) - M(t+h)]}  \\
   0  \\
\end{array}} \right).
$$
\end{Lemma}

\begin{Lemma}\label{xiao}
Write $\mathcal{{H}}(\alpha_1)$ as $K$. And denote
\begin{align}
A := &- U(t+h){\left( {(t+h)^2K} \right)^{\frac{s}{2}}}\left[ {V,U(-t-h)} \right]\nonumber\\
 &-\left[ {V,U(t+h)} \right]{\left( {(t+h)^2K} \right)^{\frac{s}{2}}}U(-t-h).\label{pxvdf4}
\end{align}
Then there holds
$$\left[ {i{\partial _t}{E_2} - K,{{\left| {{J_V}} \right|}^s}} \right] = \frac{{is}}{t+h}{\left| {{J_V}} \right|^s} + \frac{i}{t+h}U(t+h)\left[ {x \cdot \nabla ,{{\left( {(t+h)^2K} \right)}^{\frac{s}{2}}}} \right]U(-t-h) + A.
 $$
\end{Lemma}
\begin{proof}
Since $K$ is independent of $t$, Lemma \ref{50} yields
\begin{align*}
&\left[ {i{\partial _t}{E_2} - K,{{\left| {{J_V}} \right|}^s}} \right]= \left[ {i{\partial _t}{E_2} - K,U(t+h)} \right]{\left( {(t+h)^2K} \right)^{\frac{s}{2}}}U(-t-h)\\
&+ U(t+h)\left[ {i{\partial _t}{E_2} - K,{{\left( {(t+h)^2K} \right)}^{\frac{s}{2}}}U(-t-h)} \right] \\
&= \frac{i}{{2(t+h)}}{\left| {{J_V}} \right|^s} + U(t+h){\left( {(t+h)^2K} \right)^{\frac{s}{2}}}\left[ {i{\partial _t}{E_2} - K,U(-t-h)} \right]\\
&- \left[ {V,U(t+h)} \right]{\left( {(t+h)^2K} \right)^{\frac{s}{2}}}U(-t-h)+U(t+h)[i\partial_tE_2-K,(t^2K)^{\frac{s}{2}}]U(-t-h)\\
&+ \Psi{\left( {(t+h)^2K} \right)^{\frac{s}{2}}}U(-t-h),
 \end{align*}
 where the matrix $\Psi$ is
 {\small\begin{align*}
\Psi:=\left( \begin{array}{l}
 M(t+h)\frac{{ix \cdot \nabla }}{t+h} \\
  \\
 \end{array} \right.\left. \begin{array}{l}
  \\
 M(-t-h)\frac{{ix \cdot \nabla }}{t+h} \\
 \end{array} \right).
  \end{align*}}
Furthermore, by Lemma \ref{huziji9}, Lemma \ref{51} and (\ref{pxvdf4}), we get
 \begin{align*}
 \left[ {i{\partial _t}{E_2} - K,{{\left| {{J_V}} \right|}^s}} \right]&=\frac{i}{{2(t+h)}}{\left| {{J_V}} \right|^s}+\frac{{is}}{t+h}{\left| {{J_V}} \right|^s}+A+ U(t+h){\left( {(t+h)^2K} \right)^{\frac{s}{2}}}\mathfrak{M}\\
 &+ \Psi{\left( {(t+h)^2K} \right)^{\frac{s}{2}}}U(-t-h).
 \end{align*}
 where the matrix $\mathfrak{M}$ is given by  Lemma \ref{huziji9}.
Part of the terms in above formula can be cancellated. Thus we obtain
 \begin{align*}
 \left[ {i{\partial _t}{E_2} - K,{{\left| {{J_V}} \right|}^s}} \right]&= \frac{{is}}{t+h}{\left| {{J_V}} \right|^s} +A+ \Psi{\left( {(t+h)^2K} \right)^{\frac{s}{2}}}U(-t-h) - U(t+h){\left( {(t+h)^2K} \right)^{\frac{s}{2}}}\Phi\\
 &+ U(t+h){\left( {(t+h)^2K} \right)^{\frac{s}{2}}} \Upsilon\\
 &= \frac{{is}}{t+h}{\left| {{J_V}} \right|^s} + A+ \frac{i}{t+h}U(t+h)\left[ {x \cdot \nabla ,{{\left( {(t+h)^2K} \right)}^{\frac{s}{2}}}U(-t-h)} \right]\\
 &+ U(t+h){\left( {(t+h)^2K} \right)^{\frac{s}{2}}}\Upsilon,
\end{align*}
where we denote
{\small\begin{align*}
\Phi:=\left( \begin{array}{l}
 M(-t-h)\frac{{ix \cdot \nabla }}{t+h} \\
  \\
 \end{array} \right.\left. \begin{array}{l}
  \\
 M(t+h)\frac{{ix \cdot \nabla }}{t+h} \\
 \end{array} \right),
\Upsilon:=\left( \begin{array}{l}
 M(- t-h)( - \frac{{{x^2}}}{{2{(t+h)^2}}}) \\
  \\
 \end{array} \right.\left. \begin{array}{l}
  \\
 M(t+h)\frac{{{x^2}}}{{2{(t+h)^2}}} \\
 \end{array} \right)
\end{align*}}
Then our lemma follows by
\begin{align*}
 &\left[ {x \cdot \nabla ,{{\left( {(t+h)^2K} \right)}^{\frac{s}{2}}}U(-t-h)} \right]\\
 &= \left[ {x \cdot \nabla ,{{\left( {(t+h)^2K} \right)}^{\frac{s}{2}}}} \right]U(-t-h)+ {\left( {(t+h)^2K} \right)^{\frac{s}{2}}}\left[ {x \cdot \nabla ,U(-t-h)} \right] \\
 &= \left[ {x \cdot \nabla ,{{\left( {(t+h)^2K} \right)}^{\frac{s}{2}}}} \right]U(-t-h)- {\left( {(t+h)^2K} \right)^{\frac{s}{2}}}\Upsilon.
 \end{align*}
\end{proof}

\noindent Denote
\begin{align}\label{kfgtr56k}
B(s) = s{\left( {K} \right)^{\frac{s}{2}}} + \left[ {x \cdot \nabla ,{{\left( {K} \right)}^{\frac{s}{2}}}} \right],
\end{align}
then Lemma \ref{xiao} yields
\begin{align}\label{kk}
\left[ {i{\partial _t}{E_2} - K,{{\left| {{J_V}} \right|}^s}} \right] = i(t+h)^{s-1}U(t+h)B(s)U(-t-h) + A.
\end{align}

\noindent For $B(s)$, we have
\begin{Lemma}\label{55}
Let $r=1^+$, for any $g\in L^2$, we have ${\left\| {B(s)\mathbf{g}} \right\|_{r}} \le C{\left\| \mathbf{g }\right\|_2 }.$
\end{Lemma}
\begin{proof}
A small modification of the proof of Lemma 7.1 in \cite{CGV} leads to
$$B(s)=c\int_{\Gamma}\tau^{\frac{s}{2}}(\tau-K)^{-1}V_3(\tau-K)^{-1}d\tau.$$
where $V_3=2V+x\frac{d}{dx}V$, $K=\mathcal{H}(\alpha_1)$. The detailed proof is given in Appendix A.
Similar arguments as Lemma \ref{22} with additional efforts yield our lemma.
In fact, as presented in the proof of Lemma \ref{22}, we split $\Gamma$ into four parts.
Define
$${B_j}\mathbf{g} = \int_{{\Gamma _j}} {{\lambda ^{\frac{s}{2}}}{{\left( {\lambda  - K} \right)}^{ - 1}}} {V_3}{\left( {\lambda  - K} \right)^{ - 1}}\mathbf{g}d\lambda. $$
Let $\frac{1}{r}=\frac{1}{2}+\frac{1}{q}$. For $B_1$,  Lemma \ref{41} and H\"older inequality give
\begin{align*}
{\left\| {{B_1}\mathbf{g}} \right\|_r} &\le {\left\| {\int_{{\Gamma _1}} {{\lambda ^{\frac{s}{2}}}{{\left( {\lambda  - K} \right)}^{ - 1}}} {V_3}{{\left( {\lambda  - K} \right)}^{ - 1}}\mathbf{g}d\lambda } \right\|_r} \\
& \le \int_{{\Gamma _1}} {{{\left\| {{{\left( {\lambda  - K} \right)}^{ - 1}}} \right\|}_{r \to r}}{\lambda ^{\frac{s}{2}}}} {\left\| {{V_3}} \right\|_q}{\left\| {{{\left( {\lambda  - K} \right)}^{ - 1}}} \right\|_{2 \to 2}}d\lambda {\left\| \mathbf{g} \right\|_2} \\
&\le C{\left\| \mathbf{g} \right\|_2}\int_{{\Gamma _1}} {{\lambda ^{\frac{s}{2} - 2}}} d\lambda
 \le C{\left\| \mathbf{g }\right\|_2}.
\end{align*}
In the following, all $L^p$ norms refer to $L^p_{\rm{rad}}$, i.e. we only consider radial functions. Denote the resolvent set of $K$ in $L^2_{rad}$ by ${\rho _{{L^2}}}(K)$.  The resolvent set, the spectrum of $K$ in ${{L^2} \cap {L^r}}$ are denoted by ${\rho _{{L^2} \cap {L^r}}}(K)$ and ${\sigma _{{L^2} \cap {L^r}}}(K)$ respectively.
The domain of $K$ in ${{L^2} \cap {L^r}}$ is taken as $W^{2,r}\cap W^{2,2}$. We claim ${\rho _{{L^2}}}(K) \subseteq {\rho _{{L^2} \cap {L^r}}}(K)$.
In fact, let $\mathbf{g}\in L^r\cap L^2$, $\lambda\in {\rho _{{L^2}}}(K)$, the definition of resolvent set indicates there exists $\mathbf{f}\in L^2$ such that $\mathbf{g} = \left( {K - \lambda } \right)\mathbf{f}$, i.e.,
\begin{align}\label{zxc}
\left\{ \begin{array}{l}
  - \Delta {f_1} + {\tau _1}{f_1} - \lambda {f_1} + {V_1}({\alpha _1}){f_1} - {V_2}({\alpha _1}){f_2} = {g_1} \\
 \Delta {f_2} - {\tau _1}{f_2} - \lambda {f_2} + {V_2}({\alpha _1}){f_1} - {V_1}({\alpha _1}){f_2} = {g_2} \\
 \end{array} \right.
\end{align}
H\"older inequality implies
${g_1}-{V_1}({\alpha _1}){f_1} + {V_2}({\alpha _1}){f_2}\in L^r$. Therefore we have
$$- \Delta {f_1} + ({\tau _1} - \lambda) {f_1}\in L^r.$$
Since ${\tau _1} - \lambda \notin (-\infty,0]$, $- \Delta  + ({\tau _1} - \lambda)$ is invertible in $L^r$, which can be directly proved by checking the integral kernel of the resolvent of $\Delta$. Hence $f_1\in L^r$, again by the first equation in (\ref{zxc}) we have $\Delta {f_1}\in L^r$, thus $f_1\in W^{2,r}$. The same arguments show $f_2\in W^{2,r}$, by which we have proved our claim. Thus we have
${\sigma _{{L^2} \cap {L^r}}}(K) \subseteq {\sigma _{{L^2}}}(K)$, which implies for $\lambda\in {\Gamma _2} \cup {\Gamma _3} \cup {\Gamma _4}$,
$$
{\left\| {{{\left( {\lambda  - K} \right)}^{ - 1}}} \right\|_{{L^r} \cap {L^2} \to {L^r} \cap {L^2}}} \le C.
$$
Hence by H\"older inequality and the boundedness of $\Gamma_j$ where $j\in \{2,3,4\}$, we obtain for $\frac{1}{q}+\frac{1}{2}=\frac{1}{r}$
\begin{align*}
{\left\| {{B_j}\mathbf{g}} \right\|_r} &\le {\left\| {\int_{{\Gamma _j}} {{\lambda ^{\frac{s}{2}}}{{\left( {\lambda  - K} \right)}^{ - 1}}} {V_3}{{\left( {\lambda  - K} \right)}^{ - 1}}\mathbf{g}d\lambda } \right\|_r} \\
&\le \int_{{\Gamma _j}} {{\lambda ^{\frac{s}{2}}}{{\left\| {{{\left( {\lambda  - K} \right)}^{ - 1}}} \right\|}_{{L^r} \cap {L^2} \to {L^r} \cap {L^2}}}} \left\| V_3 \right\|_{L^q\bigcap L^{\infty}}{\left\| {{{\left( {\lambda  - K} \right)}^{ - 1}}} \right\|_{2 \to 2}}d\lambda {\left\| \mathbf{g} \right\|_2} \\
&\le C{\left\| \mathbf{g} \right\|_2}.
\end{align*}
\end{proof}

\begin{Lemma}\label{hua}
Recall $K=\mathcal{H}(\alpha_1)$ and $K^\frac{s}{2}$ defined via  (2.2). Let $r=1^+$, then for $\forall g\in L^2$, we have ${\left\| {{K^{\frac{s}{2}}}\mathbf{g} - \mathcal{H}_0^{\frac{s}{2}}\mathbf{g}} \right\|_r} \le C{\left\| \mathbf{g }\right\|_2}.$
\end{Lemma}
\begin{proof}
From the proof of Lemma \ref{22},
$$
{\left\| {{K^{\frac{s}{2}}}\mathbf{g} - \mathcal{H}_0^{\frac{s}{2}}\mathbf{g}} \right\|_r} = {\left\| {\int_\Gamma  {{\lambda ^{\frac{s}{2}}}{{\left( {\lambda  - K} \right)}^{ - 1}}V{{\left( {\lambda  - {{H}_0}} \right)}^{ - 1}}\mathbf{g}} } d\lambda\right\|_r}.
$$
The same arguments in the proof of Lemma \ref{55} yield our Lemma.
\end{proof}

\section{Close the bootstrap via Nonlinear estimates }

Let $0<\epsilon\ll h^{-4}\ll 1.$ $h$ is determined first, $\epsilon$ is then determined according to $h$.
Define ${\left\| u \right\|_\Sigma } = {\left\| u \right\|_{{H^2}}} + {\left\| {\left| x \right|{\partial _x}u} \right\|_2} + {\| {{{\| x |}^2}u}\|_2}$. And assume $\|g(0,x)\|_{\Sigma}\le \epsilon.$

\begin{Lemma}\label{56}
Recall $f$ and $g$ defined by (\ref{777}), (\ref{kiou8}), and ${\bf g}=(g,\bar{g})^t$. Recall the definition of $A$ in (\ref{pxvdf4}). Let $r=1^+$, and $(p',r')$ be a Strichartz admissible pair. Then for $t_1\in\mathcal{A}$,
$${\left\| {A\mathbf{g}} \right\|_{L_t^{{p}}[0,t_1] L_x^{r}}} \le C(h^{-1}){\left\| {{P_2}{{\left| {{J_V}(\alpha_1)} \right|}^s}\mathbf{g}} \right\|_{L_t^\infty[0,t_1]L_x^2}}.$$
\end{Lemma}
\begin{proof}
Denote
$$\widetilde{V}(t) = [M(t+h) - M( -t-h)]\left( {\begin{array}{*{20}{c}}
   0  \\
   {{V_2}}  \\
\end{array}} \right.\left. {\begin{array}{*{20}{c}}
   {-{V_2}}  \\
   0  \\
\end{array}} \right),  \widehat{V} (t) = \left( {\begin{array}{*{20}{c}}
   0  \\
   {{V_2}}  \\
\end{array}} \right.\left. {\begin{array}{*{20}{c}}
   {-{V_2}}  \\
   0  \\
\end{array}} \right).
 $$
Let $\frac{1}{r}=\frac{1}{2}+\frac{1}{q}$.
Lemma \ref{hua}, fractional Leibnitz formula (Appendix of \cite{KPV}), Lemma \ref{22}, Lemma \ref{31} and Lemma \ref{a2} imply
{\small\begin{align*}
{\left\| {A\mathbf{g}} \right\|_{r}}&\le \left\| {(t+h)^s}\mathcal{H}_0^{\frac{s}{2}}[U(t+h)\widetilde{V}(t)U( -t-h)\mathbf{g}] \right\|_{r} + \left\| \widetilde{V} (t){{\left| {{J_V}} \right|}^s}\mathbf{g} \right\|_{r} + {(t+h)^s}\left\| \widetilde{V}(t)\mathbf{g} \right\|_2 \\
&\le {\left\| {(t+h)^sU(t+h)\widetilde{V}(t)U( - t-h)\mathbf{g}} \right\|_{{W^{s,r}}}} + \left\| \widehat{V}(t)\left| x \right|^2 \right\|_q\left\| \left| {{J_V}} \right|^s\mathbf{g} \right\|_2{(t+h)^{ - 1}}\\
&+ {(t+h)^{s-1}}{\left\|\mathbf{g} \right\|_2}\left\| \widehat{V} (t)\left| x \right|^2 \right\|_{\infty}\\
&\lesssim {\left\| {{(t+h)^s}U(t+h)\widetilde{V}(t)} \right\|_{{W^{s,2}}}}{\left\| \mathbf{g} \right\|_q} + {\left\| {{(t+h)^s}U( - t-h)\mathbf{g}} \right\|_{{W^{s,2}}}}\left\|
\widetilde{V} (t) \right\|_q\\
&+ {\left\| {{{\left| {{J_V(\alpha_1)}} \right|}^s}\mathbf{g}} \right\|_2}{(t+h)^{ - 1}} + {(t+h)^{s-1}}{\left\| \mathbf{g} \right\|_2}\\
&\lesssim (t+h)^{s-1}{\left\| \mathbf{g} \right\|_q} +(t+h)^{-1}{\left\| {{{\left| {{J_V}(\alpha_1)} \right|}^s}\mathbf{g}} \right\|_2} + (t+h)^{s-1}{\left\| \mathbf{g} \right\|_2} \\
&\lesssim  (t+h)^{-1}{\left\| {{{\left| {{J_V}(\alpha_1)} \right|}^s}\mathbf{g}} \right\|_2} \le C(t+h)^{-1}{\left\| {{P_2}{{\left| {{J_V(\alpha_1)}} \right|}^s}\mathbf{g}} \right\|_{L_t^\infty L_x^2}}.
\end{align*}}
\end{proof}

\noindent Recall $B(s)$ defined in (\ref{kfgtr56k}). Denote $\bar{B}(s)\mathbf{g}={i{(t+h)^{s-1}}U(t+h)B(s)U( -t-h)\mathbf{g}}$.
\begin{Lemma}\label{57}
Let $r=1^+$, and $(p',r')$ be an admissible pair, then for $t_1\in\mathcal{A}$
$${\left\| \bar{B}(s)\mathbf{g} \right\|_{L_t^{p}[0,t_1]L_x^{r}}} \le C(h^{-1}){\left\| {{P_2}{{\left| {{J_V(\alpha_1)}} \right|}^s}\mathbf{g}} \right\|_{L_t^\infty[0,t_1] L_x^2}}.$$
\end{Lemma}
\begin{proof}
Thanks to Lemma \ref{55}, Lemma \ref{31} and Lemma \ref{a2}, it is easy to see
\begin{align*}
&{\left\| {i(t+h)^{s-1}U(t+h)B(s)U( -t-h)\mathbf{g}} \right\|_{L_t^{p}L_x^{r}}}\\
&\le {\left\| {i{(t+h)^{s-1}}\mathbf{g}} \right\|_{L_t^{p}L_x^2}} \le \left\| {(t+h)^{s-1}}\right\|_{L_t^{r}}{\left\| {{{\left| {{J_V}(\alpha_1)} \right|}^s}\mathbf{g}} \right\|_{L_t^\infty L_x^2}}\\
&\le C(h^{-1}){\left\| {{P_2}{{\left| {{J_V}(\alpha_1)} \right|}^s}h(t)} \right\|_{L_t^\infty L_x^2}}.
\end{align*}
\end{proof}

Let us estimate the $\{{\bf D}_i\}^{4}_{i=0}$ terms in (\ref{b1}). The main tools are  Lemma \ref{hua}, Lemma \ref{31}, and fractional Leibnitz formula (Appendix of \cite{KPV}).
\begin{Lemma}\label{80}
Let $r=1^+$, and $(p',r')$ be an admissible pair. Then for $t_1\in\mathcal{A}$, we have (all the spacetime norms used below are restricted in $[0,t_1]\times\Bbb R$)
{\small \begin{align*}
&{\left\| {{{\left| {{J_V}(\alpha_1)} \right|}^s}\bf{D}_0} \right\|_{L_t^{p}L_x^{r}}} \lesssim \left\| {{P_2}{{\left| {{J_V}(\alpha_1)} \right|}^s}\mathbf{g}} \right\|_{L_t^\infty L_x^2}^2\\
&{\left\| {{{\left| {{J_V}(\alpha_1)} \right|}^s}\bf{D}_4} \right\|_{L_t^{p}L_x^{r}}} \lesssim \left\|P_2 {{{\left| {{J_V}(\alpha_1)} \right|}^s}\mathbf{g}} \right\|_{L_t^\infty L_x^2}^2 + \left\| P_2{{{\left| {{J_V}(\alpha_1)} \right|}^s}\mathbf{g}} \right\|_{L_t^\infty L_x^2}^{m - 1} + \left\| P_2{{{\left| {{J_V}(\alpha_1)} \right|}^s}\mathbf{g}} \right\|_{L_t^\infty L_x^2}^{m - 2}\\
&{\left\| {{{\left| {{J_V}(\alpha_1)} \right|}^s}\bf{D}_j} \right\|_{L_t^{p}L_x^{r}}} \lesssim {\left\| P_2{{{\left| {{J_V}(\alpha_1)} \right|}^s} \mathbf{g}} \right\|_{L_t^\infty L_x^2}^3},
\end{align*} }
for $j=1,2,3.$
\end{Lemma}
\begin{proof} By (\ref{zx}),
\begin{align*}
&{\left\| {{{\left| {{J_V}(\alpha_1)} \right|}^s}\bf{D}_0} \right\|_{L_t^pL_x^{r}}} \le C{\left\| {\left| {\gamma '} \right|{(t+h)^s} + \left| {\omega '} \right|{(t+h)^s}} \right\|_{L_t^{p}}} \le C\mathcal{M}^2{\left\| {{t^{ - 1}}} \right\|_{L_t^{p}}} \le C\mathcal{M}^2.
\end{align*}
We can write $D_4$ as $D_4=e^{-i\Omega}\int_0^1 [ {A_1} + {A_2} + {A_3}+{A_4}+{A_5}](1 - \tau)d\tau$, where
{\small \begin{align*}
{A_1} &= {\tilde g ^2}{\left( {\overline {\varphi  + \tau \tilde g} } \right)^2}(\varphi  + \tau \tilde g)F''({\left| {\varphi  + \tau \tilde g} \right|^2}), \mbox{ }
 {A_2}  = 2{\left| \tilde g  \right|^2}{(\varphi  + \tau \tilde g)^2}\overline {\varphi  + \tau \tilde g} F''({\left| {\varphi  + \tau \tilde g} \right|^2}), \\
 {A_3} &= {\left( {\bar {\tilde g} } \right)^2}{(\varphi  + \tau \tilde g)^3}F''({\left| {\varphi  + \tau \tilde g} \right|^2}), \mbox{ }
 {A_4}  = 2{\tilde g ^2}\overline {\varphi  + \tau \tilde g} F'({\left| {\varphi  + \tau \tilde g} \right|^2}), \\
 {A_5} &= 4{\left| \tilde g  \right|^2}(\varphi  + \tau \tilde g)F'({\left| {\varphi  + \tau \tilde g} \right|^2}), \mbox{ }
\tilde g = g{e^{i\Omega }}.
\end{align*}}
We give the detailed proof for $A_3$, the proof of $A_1,A_2,A_4,A_5$ are almost the same. Let $\frac{1}{r}=\frac{1}{2}+\frac{1}{q}$, $r=1^+$, $\left\{ \eta  \right\}\triangleq {\left( {\eta ,\bar \eta } \right)^t}$. By Lemma \ref{hua},
{\footnotesize \begin{align*}
&{\left\| {{{\left| {{J_V}(\alpha_1)} \right|}^s}\left\{ {{{\left( {\bar{\tilde g}} \right)}^2}{{(\varphi  + \tau \tilde g)}^3}F''({{\left| {\varphi  + \tau \tilde g} \right|}^2})} \right\}} \right\|_r} \\
&\lesssim {(t+h)^s}{\left\| {U( - t-h)\left\{ {{{\left( {\bar{\tilde g}} \right)}^2}{{(\varphi  + \tau \tilde g)}^3}F''({{\left| {\varphi  + \tau \tilde g} \right|}^2})} \right\}} \right\|_{{W^{s,r}}}}\\
&+ {(t+h)^s}{\left\| {\left\{ {{{\left( {\bar{\tilde g}} \right)}^2}{{(\varphi  + \tau \tilde g)}^3}F''({{\left| {\varphi  + \tau \tilde g} \right|}^2})} \right\}} \right\|_2}\\
&\lesssim{(t+h)^s}{\left\| {M( - 3(t+h)){{(\varphi  + \tau \tilde g)}^3}F''({{\left| {\varphi  + \tau \tilde g} \right|}^2})(M(t+h)\bar g)(M(t+h)\bar g)} \right\|_{{W^{s,r}}}}\\
&+ {(t+h)^s}\left\| g \right\|_\infty ^2.
\end{align*}}
Then by fractional Leibnitz formula (Appendix of \cite{KPV}), Lemma \ref{hua}, we obtain
{\footnotesize \begin{align*}
&{\left\| {{{\left| {{J_V}} \right|}^s}\left\{ {{{\left( {\bar{\tilde g}} \right)}^2}{{(\varphi  + \tau \tilde g)}^3}F''({{\left| {\varphi  + \tau \tilde g} \right|}^2})} \right\}} \right\|_r}\\
&\lesssim {(t+h)^s}{\left\| {M( - 3(t+h)){{(\varphi  + \tau \tilde g)}^3}F''({{\left| {\varphi  + \tau \tilde g} \right|}^2})} \right\|_{{H^2}}}{\left\| {{g^2}} \right\|_q} + {(t+h)^s}\left\| g \right\|_\infty ^2 \\
&+ {(t+h)^s}{\left\| {M( - 3(t+h)){{(\varphi  + \tau \tilde g)}^3}F''({{\left| {\varphi  + \tau \tilde g} \right|}^2})} \right\|_q}{\left\| {{{(M(t+h)\bar g)}^2}} \right\|_{{H^s}}}\\
&\lesssim {(t+h)^s}{\left\| {M( - 3(t+h)){{(\varphi  + \tau \tilde g)}^3}F''({{\left| {\varphi  + \tau \tilde g} \right|}^2})} \right\|_{{H^2}}}{\left\| g \right\|_\infty }{\left\| g \right\|_q}\\
&+ {(t+h)^s}{\left\| {M(t+h)\bar g} \right\|_{{H^s}}}{\left\| g \right\|_\infty }+ {(t+h)^s}\left\| g \right\|_\infty ^2 \\
&\lesssim {(t+h)^s}{\left\| {M( - 3(t+h)){{(\varphi  + \tau \tilde g)}^3}F''({{\left| {\varphi  + \tau \tilde g} \right|}^2})} \right\|_{{H^2}}}{\left\| g \right\|_\infty }{\left\| g \right\|_q}\\
&+ {\left\| {{{\left| {{J_V}} \right|}^s}\mathbf{g}} \right\|_2}{\left\| g \right\|_\infty } +{(t+h)^s}\left\| g \right\|_\infty ^2.
\end{align*} }
Thus Lemma \ref{31}  gives
 {\footnotesize\begin{align}
& {\left\| {{{\left| {{J_V}} \right|}^s}\left\{ {{{\left( {\bar{\tilde g}} \right)}^2}{{(\varphi  + \tau \tilde g)}^3}F''({{\left| {\varphi  + \tau \tilde g} \right|}^2})} \right\}} \right\|_r}\label{Ghyu7nnjk}\\
&\lesssim {(t+h)^{s - 2}}{\left\| {{{\left| x \right|}^2}{{(\varphi  + \tau \tilde g)}^3}F''({{\left| {\varphi  + \tau \tilde g} \right|}^2})} \right\|_2}{\left\| g \right\|_\infty }{\left\| g \right\|_q}\nonumber\\
&+{(t+h)^{s-1}}{\left\| {\left| x \right|{{\left( {{{(\varphi  + \tau \tilde g)}^3}F''({{\left| {\varphi  + \tau \tilde g} \right|}^2})} \right)}^\prime }} \right\|_2}{\left\| g \right\|_\infty }{\left\| g \right\|_q} \nonumber\\
&\mbox{ }\mbox{ }+ {(t+h)^s}{\left\| {{{\left( {{{(\varphi  + \tau \tilde g)}^3}F''({{\left| {\varphi  + \tau \tilde g} \right|}^2})} \right)}^{\prime \prime }}} \right\|_2}{\left\| g \right\|_\infty }{\left\| g \right\|_q}\nonumber\\
&+ {\left\| {{{\left| {{J_V}} \right|}^s}\mathbf{g}} \right\|_2}{\left\| g \right\|_\infty } + {(t+h)^s}\left\| g \right\|_\infty ^2.\nonumber
\end{align} }
Then the assumption of $F$ shows the LHS of (\ref{Ghyu7nnjk}) is bounded by
{\footnotesize \begin{align*}
& {(t+h)^{s - 2}}\left( {{{\left\| g \right\|}_\infty }{{\left\| g \right\|}_q} + {{\left\| {{{\left| x \right|}^2}g} \right\|}_2}{{\left\| g \right\|}_q}\left( {\left\| g \right\|_\infty ^{m - 2} + \left\| g \right\|_\infty ^{n - 2}} \right)} \right)\\
&+ {(t+h)^{s - 1}}\left( {{{\left\| g \right\|}_\infty }{{\left\| g \right\|}_q}+ {{\left\| {\left| x \right|g} \right\|}_2}{{\left\| g \right\|}_q}\left( {\left\| g \right\|_\infty ^{m - 3}+ \left\| g \right\|_\infty ^{n - 3}} \right)} \right) \\
&+ {\left\| {{{\left| {{J_V}} \right|}^s}\mathbf{g}} \right\|_2}{\left\| g \right\|_\infty }+ {(t+h)^s}\left\| g \right\|_\infty ^2 + {(t+h)^s}{\left\| g \right\|_\infty }{\left\| g \right\|_q}
\end{align*} }
which by  Lemma \ref{qqq} is further dominated by
{\footnotesize \begin{align*}
& {(t+h)^{ - s}}\left\| {{{\left| {{J_V}} \right|}^s}\mathbf{g}} \right\|_2^2 + {(t+h)^{s + 2 - s(m - 1)}}\left\| {{{\left| {{J_V}} \right|}^s}\mathbf{g}} \right\|_2^{m - 1}\\
&+ {(t+h)^{s - (m - 2)s}}\left\| {{{\left| {{J_V}} \right|}^s}\mathbf{g}} \right\|_2^{m - 2} + {(t+h)^{s + 2 - s(n - 1)}}\left\| {{{\left| {{J_V}} \right|}^s}\mathbf{g}} \right\|_2^{n - 1}\\
& + {(t+h)^{s - (n - 2)s}}\left\| {{{\left| {{J_V}} \right|}^s}\mathbf{g}} \right\|_2^{n - 2}.
\end{align*} }
Since $s=\frac{7}{4}^+$, $m,n>\frac{26}{7}$, we get the desired estimates of $\mathbf{D}_4$ in Lemma \ref{80}.

Since the proofs of $\mathbf{D}_1$, $\mathbf{D}_2$, $\mathbf{D}_3$  are almost the same, we only give the proof of $\mathbf{D}_3$. By the definition of $\mathbf{D}_3$,
\begin{align*}
 &{\left\| {{{\left| {{J_V}} \right|}^s}{\mathbf{D}_3}} \right\|_r}\lesssim {(t+h)^s}{\left\| {{K^{\frac{s}{2}}}U( - t-h){\mathbf{D}_3}} \right\|_r}.
 \end{align*}
And Lemma \ref{hua}, fractional Leibnitz formula(Appendix of \cite{KPV}) and Lemma \ref{31} yield
{\footnotesize\begin{align*}
&{\left\| {{K^{\frac{s}{2}}}U( - t-h){\mathbf{D}_3}} \right\|_r}\\
&\lesssim \int_{{\alpha _1}}^{\alpha (t)} {{{\left\| {U( - t - h){\rm{ }}\left\{ {[F''({\varphi ^2}(\tau )){\varphi ^3}(\tau ){\varphi _\alpha }(\tau ) + F'({\varphi ^2}(\tau ))\varphi (\tau ){\varphi _\alpha }(\tau )]\bar g} \right\}} \right\|}_{{W^{s,r}}}}d\tau } \\
 &+\int_{{\alpha _1}}^{\alpha (t)} {{{\left\| {U( - t - h){\rm{ }}\left\{ {[F''({\varphi ^2}(\tau )){\varphi ^3}(\tau ){\varphi _\alpha }(\tau ) + F'({\varphi ^2}(\tau ))\varphi (\tau ){\varphi _\alpha }(\tau )]\bar g} \right\}} \right\|}_{{L^2}}}d\tau } \\
 &\lesssim {\int_{{\alpha _1}}^{\alpha (t)} {\left\| {M( - 2(t+h))[F''({\varphi ^2}(\tau )){\varphi ^3}(\tau ){\varphi _\alpha }(\tau ) + F'({\varphi ^2}(\tau ))\varphi (\tau ){\varphi _\alpha }(\tau )]M(t+h)\bar g} \right\|} _{{W^{s,r}}}}d\tau  \\
 &+ \left| {\alpha (t) - {\alpha _1}} \right|{\left\| g \right\|_{L^2}}
  \end{align*}}
 which by H\"older inequality is  further dominated by
{\footnotesize\begin{align*}
 &{\int_{{\alpha _1}}^{\alpha (t)} {\left\| {M( - 2(t+h))[F''({\varphi ^2}(\tau )){\varphi ^3}(\tau ){\varphi _\alpha }(\tau ) + F'({\varphi ^2}(\tau ))\varphi (\tau ){\varphi _\alpha }(\tau )]} \right\|} _{{W^{s,2}}}}{\left\| g \right\|_{L^q}}d\tau  \\
 &+ {\int_{{\alpha _1}}^{\alpha (t)} {\left\| {M( - 2(t+h))[F''({\varphi ^2}(\tau )){\varphi ^3}(\tau ){\varphi _\alpha }(\tau ) + F'({\varphi ^2}(\tau ))\varphi (\tau ){\varphi _\alpha }(\tau )]} \right\|} _{L^q}}{\left\| {M(t+h)\bar g} \right\|_{{W^{s,2}}}}d\tau  \\
&+ \left| {\alpha (t) - {\alpha _1}} \right|{\left\| g \right\|_{L^2}}.
\end{align*}}
Thus we conclude by Lemma \ref{31}, Lemma \ref{a2},
{\footnotesize\begin{align*}
{\left\| {{{\left| {{J_V}} \right|}^s}{\mathbf{D}_3}} \right\|_r}&\lesssim {(t+h)^s}\left| {\alpha (t) - {\alpha _1}} \right|{\left\| g \right\|_q} + \left| {\alpha (t) - {\alpha _1}} \right|{\left\| {{{\left| {{J_V}} \right|}^s}\mathbf{g}} \right\|_2} + {(t+h)^s}\left| {\alpha (t) - {\alpha _1}} \right|{\left\| g \right\|_2} \\
&\lesssim \left| {\alpha (t) - {\alpha _1}} \right|{\left\| {{{\left| {{J_V}} \right|}^s}\mathbf{g}} \right\|_2} \lesssim {t^{ - s + 1}}{\left\| {{{\left| {{J_V}} \right|}^s}\mathbf{g}} \right\|^3_{L_t^\infty L_x^2}}\lesssim {t^{ - s + 1}}{\left\| {{{P_2\left| {{J_V}} \right|}^s}\mathbf{g}} \right\|^3_{L_t^\infty L_x^2}},
\end{align*}}
where we have used (\ref{s1}) in the last inequality.
\end{proof}

\subsection{Proof of main Theorem}

Recall that ${\left\| u \right\|_\Sigma } = {\left\| u \right\|_{{H^2}}} + {\left\| {\left| x \right|{\partial _x}u} \right\|_2} +\| |x |^2u \|_2$, and
$P_2=P_c(H(\alpha_1))$, $K=H(\alpha_1)$.
Applying $|J_V(\alpha_1)|^s$ and $P_2$ to (\ref{b1}), then by the commutator relation (\ref{kk}), we obtain
$$
i{\partial _t}{P_2}{\left| {{J_V(\alpha_1)}} \right|^s}\mathbf{g} - K{P_2}{\left| {{J_V(\alpha_1)}} \right|^s}\mathbf{g} - {P_2}{\left| {{J_V(\alpha_1)}} \right|^s}\mathbf{D }= {P_2}\left( {\bar B(s)\mathbf{g} + A \mathbf{g}} \right).
$$
Let $(p',r')$ be an admissible pair with $4<p\le\infty$, $r=1^+$.
Let $t_1\in\mathcal{A}$. By Lemma \ref{56}, Lemma \ref{57}, Lemma \ref{80}, and Strichartz estimates (see Corollary 7.3 of Krieger, Schlag \cite{KS}), we have
{\footnotesize\begin{align*}
\left\| {{P_2}{{\left| {{J_V(\alpha_1)}} \right|}^s}\mathbf{g}} \right\|_{L_t^\infty[0,t_1] L_x^2}&\lesssim {\left\| {{{\left| {{J_V(\alpha_1)}} \right|}^s}\mathbf{D}} \right\|_{L_t^p[0,t_1]L_x^r}} + {\left\| {\bar B(s)\mathbf{g} + A\mathbf{g}} \right\|_{L_t^p[0,t_1]L_x^r}}+h^s\|g(0,x)\|_{\Sigma} \\
&\lesssim\left\| {{P_2}{{\left| {{J_V(\alpha_1)}} \right|}^s}\mathbf{g}} \right\|_{L_t^\infty [0,t_1]L_x^2}^3 + \left\| {{P_2}{{\left| {{J_V(\alpha_1)}} \right|}^s}\mathbf{g}} \right\|_{L_t^\infty[0,t_1] L_x^2}^2 \\
&+ \left\| {{P_2}{{\left| {{J_V(\alpha_1)}} \right|}^s}\mathbf{g}} \right\|_{L_t^\infty[0,t_1] L_x^2}^{m - 1} + \left\| {{P_2}{{\left| {{J_V(\alpha_1)}} \right|}^s}\mathbf{g}} \right\|_{L_t^\infty[0,t_1] L_x^2}^{m - 2} \\
&+ \left\| {{P_2}{{\left| {{J_V(\alpha_1)}} \right|}^s}\mathbf{g}} \right\|_{L_t^\infty[0,t_1] L_x^2}^{n - 1} + \left\| {{P_2}{{\left| {{J_V(\alpha_1)}} \right|}^s}\mathbf{g}} \right\|_{L_t^\infty[0,t_1] L_x^2}^{n - 2} \\
&+ C({h^{ - 1}}){\left\| {{P_2}{{\left| {{J_V(\alpha_1)}} \right|}^s}\mathbf{g}} \right\|_{L_t^\infty [0,t_1]{L^2_x}}} + h^s{\left\| {{g(0,x)}} \right\|_{{\Sigma}}}.
\end{align*}}
Thus choosing $h$ to be sufficiently large and $\|g(0)\|_{\Sigma}$ sufficiently small, we obtain for any $t_1\in\mathcal{A}$
\begin{align}\label{dfre5t}
\mathcal{M}_{t_1}\lesssim {\left\| {{P_2}\left| {{J_V(\alpha_1)}} \right|g} \right\|_{L_t^\infty[0,t_1] L_x^2}}<C\epsilon+C\kappa^{2}.
\end{align}
Since $\epsilon\ll \kappa$ and $\kappa\ll1$, (\ref{dfre5t}) shows $\mathcal{A}$ is open. Then for $\mathcal{A}$ is defined to be closed, we have $\mathcal{A}=[0,\infty)$ and thus for all $t\in(0,\infty)$,
\begin{align}\label{dfre5t2}
\mathcal{M}_{t}\le\kappa.
\end{align}
It is standard to deduce our main theorem from (\ref{dfre5t2}) (e.g. \cite{BP}). In fact,
(\ref{zx}) implies that $\gamma(t)$, $\omega(t)$ have limits $\gamma_{\infty}$, $\omega_{\infty}$ as $t\to\infty$. Consequently, we can introduce the limit trajectory $\sigma_+(t)$:
$$
\beta_+(t)=\omega_+t+\gamma_+, \omega_+=\omega_{\infty}, \gamma_+=\int_0^{\infty}(\omega(\tau)-\omega_{\infty})d\tau.
$$
Obviously, $\sigma(t)-\sigma_+(t)=O(t^{-s+1})$ as $t\to\infty$. Hence we obtain the limit soliton $w(x;\sigma_+(t))$ and
$$
\|w(x;\sigma(t))-w(x;\sigma_+(t))\|_{L^2\bigcap L^{\infty}}=O(t^{-s+1}).
$$
Introduce the transformation
$\chi=e^{i\Phi_{\infty}}g(x,t)$, $\Phi_{\infty}=-\beta_+(t)$.
Repeat the construction in Section 2 to Section 5 with the operator $K$, $P_2$ replaced by $K_+$ and $P_c(K_+)$ respectively, where
\begin{align*}
{K_ + } = \left( {\begin{array}{*{20}{c}}
   { - \Delta  - {\omega _ + }}  \\
   {}  \\
\end{array}} \right.\left. {\begin{array}{*{20}{c}}
   {}  \\
   {\Delta  + {\omega _ + }}  \\
\end{array}} \right) + \left( {\begin{array}{*{20}{c}}
   {{V_1}({\alpha _ + })}  \\
   {{V_2}({\alpha _ + })}  \\
\end{array}} \right.\left. {\begin{array}{*{20}{c}}
   { - {V_2}({\alpha _ + })}  \\
   { - {V_1}({\alpha _ + })}  \\
\end{array}} \right).
\end{align*}
Then we can also prove $\|\chi\|_{\infty}\le C t^{-s}$. Therefore, we have obtained
$$u=w(x;\sigma_+(t))+\chi+O(t^{-s+1}),$$
which combined with the estimate of $\chi$ yields Theorem 1.1.

\section{Appendix A}
\begin{Lemma}\label{qqq}
Let $F$ satisfy Assumption A. If the initial data $u_0\in H^2$ satisfies $\|u_0(1+|x|)^2\|_2<\infty$, then there exists a unique solution $u(t)$ to (\ref{1}), and
$$\left\|u(t)\right\|_{H^2}\le C(1+t), \mbox{  } \|u(t)(1+|x|)\|_2\le C(1+t),\mbox{  } \|u(t)(1+|x|^2)\|_2\le C(1+t^3).$$
\end{Lemma}
\begin{proof}
{\bf Step One. Growth of Sobolev norms}.
First of all we prove the Bourgain-Staffilani bound for the growth of $\|u\|_{H^2}$:
\begin{align}\label{pi}
\|u\|_{H^2}\le C(t+1).
\end{align}
(\ref{pi}) is indeed contained in the proof of Theorem 2.2 of \cite{Sta}. Althoguh  Staffilani considered the nonlinearities of the form $\Sigma_{\alpha_1+\alpha_2=m}u^{\alpha_1}{\bar{u}}^{\alpha_2}$ with $m\ge4$ in [Theorem 2.2,\cite{Sta}] , her proof can be extended to our case with very small modifications. For reader's convenience, we write down the proof in Lemma \ref{q12} below. The originality of the proof belongs to Staffilani \cite{Sta}.

{\bf Step Two. Growth of moments.} In the following $u'$ demotes $\partial_xu$ and $u''$ denotes  $\partial_{xx}u$.
First we prove $\|u|x|\|_2\le C(t+1)$. $ux$ satisfies
\begin{align}\label{788}
i{\partial _t}(xu) =  - \left( {\Delta u} \right)x - F({\left| u \right|^2})ux,
\end{align}
Multiplying (\ref{788}) by $\bar ux$, then taking the imaginary part, we obtain
$$\frac{1}{2}\frac{d}{{dt}}\left\| {\left| x \right|u} \right\|_2^2 =  2{\rm{Im}}\int_{\Bbb R} {\left( {\nabla u} \right)x\bar u}.$$
Hence one has
$$\frac{d}{{dt}}\left\| {\left| x \right|u} \right\|_2^2 \le C{\left\| {\left| x \right|u} \right\|_2}.$$
Gronwall's inequality yields $\|u|x|\|_2\le C(t+1)$.
Second, we prove
\begin{align}\label{pi1}
\|u'|x|\|_2\le C(1+t^2).
\end{align}
The equation for $u'x$ is
$$i{\partial _t}(xu') =  - \left( {\Delta u'} \right)x - F'({\left| u \right|^2})(u'\bar u + \bar u'u)ux - F({\left| u \right|^2})u'x.
$$
Multiplying the above formula with $\bar u'x$, taking the imaginary part, we have
$$\frac{d}{{dt}}\left\| {\left| x \right|u'} \right\|_2^2 \lesssim {\left\| {u''} \right\|_2}{\left\| {\left| x \right|u'} \right\|_2} + \left\| {u'} \right\|_\infty \left( {\left\| u \right\|_\infty ^{m - 2} + \left\| u \right\|_\infty ^{n - 2}} \right)\left\| {\left| x \right|u'} \right\|_2\left\| {\left| x \right|u} \right\|_2.
$$
Since $\|u\|_{H^2}$ is bounded by $C(1+t)$ and $\|u\|_{H^1}\le C$, from Sobolev imbedding theorem and (\ref{pi1}) , it follows,
$$\frac{d}{{dt}}\left\| {\left| x \right|u'} \right\|_2^2 \le C{\left\| {\left| x \right|u'} \right\|_2} (1+t).$$
Thus Gronwall's inequality implies
\begin{align}\label{pi2}
{\left\| {\left| x \right|u'} \right\|_2} \le C({t^2} + 1).
\end{align}
Third, we prove $\|u|x|^2\|_2\le C(t^3+1)$.
Similar arguments as the first step, we can show give
$$\frac{d}{{dt}}\left\| {{{\left| x \right|}^2}u} \right\|_2^2 \le C\left| {\int_{\Bbb R} {u'\bar u{x^3}} } \right| \le C{\left\| {{{\left| x \right|}^2}u} \right\|_2}{\left\| {\left| x \right|u'} \right\|_2} \le C({t^2} + 1){\left\| {{{\left| x \right|}^2}u} \right\|_2}.$$
Again by Gronwall's inequality, we get the desired result.
\end{proof}

\begin{Lemma}\label{q12}
Let $F$ satisfy Assumption A. If the initial data $u_0\in H^2$ satisfies $\|u_0\|_{H^2}<\infty$, then there exists a unique solution $u(t)$ to (\ref{1}) and
\begin{align}\label{ui}
\left\|u(t)\right\|_{H^2}\le C(1+t)
\end{align}
\end{Lemma}
\begin{proof}
The originality of the proof belongs to Staffilani \cite{Sta}. The presentation here is a very small modification of  [Theorem 2.2, \cite{Sta}, Page 137-139].
Recall the norms introduced by \cite{Sta}:
{\small\begin{align*}
&\nu^s_{1}(v):=\|\partial^{s+\frac{1}{2}}_{x}v\|_{L^{\infty}_xL^{2}_T},\mbox{  }
\nu^s_{2}(v):=\|\partial^{(s-\frac{1}{2})^{-}}_{x}v\|_{L^{2}_xL^{\infty}_T},\\
&
\nu^s_{3}(v):=\|\partial^{s}_{x}v\|_{L^{\infty}_TL^{2}_x},\mbox{  }
\Omega^s_{T}(v):=\max_{i}\nu_{i}(v)+\|v\|_{L^{\infty}_TL^2_x}.
\end{align*}}
and the Banach space
\begin{align}
X^s_{T}:=\{v\in  C([0, T], H^s):\Omega^s_{T}(v)<\infty\}.
\end{align}
{\bf Step 1. Claim A.} For any $s\ge1$, there exists a time $T$ depending only on $\|u_0\|_{H^1}$ and a unique solution to (\ref{1}) such that
\begin{align}
\Omega^s_{T}(u)\le C\|u_0\|_{H^2}.
\end{align}
{\bf Proof of Claim A}
Denoting the linear Schr\"odinger group by $S(t)$, it suffices to apply contraction argument to the operator
\begin{align}
Lu(t):=S(t)u_0+i\int^{t}_0S(t-\tau)F(|u|^2)u(\tau)d\tau.
\end{align}
[Lemma 4.1 \cite{Sta}] and Minkowski inequality give
\begin{align}\label{pi4}
\Omega_{T}( Lu)\le C\|u_0\|_{H^2}+T^{\frac{1}{2}}\|F(|u|^2)u\|_{L^2_{T}H^2}.
\end{align}
Then Assumption A(ii), fractional chain rule (\cite{KPV}) imply
\begin{align*}
\|F(|u|^2)u\|_{L^2_{T}H^s_x}&\le T^{\frac{1}{2}}\|F(|u|^2)u\|_{L^{\infty}_{T}H^s_x}\le T^{\frac{1}{2}}(\|u\|^{m}_{L^{\infty}_{T}H^s_x}+ T^{\frac{1}{2}}\|u\|^{n}_{L^{\infty}_{T}H^s_x}).
\end{align*}
Thus (\ref{pi4}) yields
\begin{align}\label{pi5}
\Omega_{T}( Lu)\le C\|u_0\|_{H^2}+T(\Omega^s_{T}(u)^m+\Omega^s_{T}(u)^n).
\end{align}
Similarly one has $L$ is a contraction if $T$ is small enough. Thus Claim A follows.

{\bf Step 2. Proof of (\ref{ui}).} Calculating $\frac{d}{dt}\|u\|^2_{{\dot H}^2_x}$ gives
\begin{align}\label{pi6}
\|u(t)\|^2_{{\dot H}^2_x}\le \|u_0\|^2_{{\dot H}^2_x}+\int^T_0\int_{\Bbb R}|\partial^{2}_x(F(|u|^2)u)\partial^2_xu|dxdt.
\end{align}
Meanwhile Claim A shows $\|u\|_{L^{2}_{x}L^{\infty}_T}\le  C\|u_0\|_{H^{\frac{1}{2}+}}$, $\|\partial_xu\|_{L^{2}_{x}L^{\infty}_T}\le C\|u_0\|_{H^{\frac{3}{2}+}}$, $\|\partial^2_xu\|_{L^{\infty}_{x}L^{2}_T}\le C\|u_0\|_{H^{\frac{3}{2}}}$,
and $\|\partial_xu\|_{L^{\infty}_{x}L^{2}_T}\le C\|u_0\|_{H^{\frac{1}{2}}}$.
Hence by Assumption A(ii), one obtains
{\small\begin{align*}
&\int^T_0\int_{\Bbb R}|\partial^{2}_x(F(|u|^2)u)\partial^2_xu|dxdt\le \|\partial^2_xu\|^2_{L^{\infty}_xL^2_T}\|u\|^{2}_{L^{2}_{x}L^{\infty}_T}(\|u\|^{m-3}_{L^{\infty}_{T}L^{\infty}_x}+ \|u\|^{n-3}_{L^{\infty}_{T}L^{\infty}_x})\\
&+\|\partial^2_xu\|_{L^{\infty}_xL^2_T}\|\partial_xu\|^{2}_{L^{\infty}_{x}L^{2}_T}
\|\partial_xu\|^{2}_{L^{2}_{x}L^{\infty}_T}\|u\|^{2}_{L^{2}_{x}L^{\infty}_T}(\|u\|^{m-3}_{L^{\infty}_{T}L^{\infty}_x}+ \|u\|^{n-3}_{L^{\infty}_{T}L^{\infty}_x}))\\
&\le  (\|u_0\|^2_{H^{\frac{3}{2}}}+\|u_0\|_{H^{\frac{3}{2}}}\|u_0\|_{H^{\frac{3}{2}+}})(\|u_0\|^{m-1}_{H^1_x}+\|u_0\|^{m-1}_{H^1_x}).
\end{align*}}
By the interpolation inequality, $\|f\|_{{H}^s}\le \|f\|^{2-s}_{{ H}^1}\|f\|^{s-1}_{{H}^2}$, we have
\begin{align*}
&\int^T_0\int_{\Bbb R}|\partial^{2}_x(F(|u|^2)u)\partial^2_xu|dxdt\le C(\|u_0\|_{H^1_x})\|u_0\|_{H^{2}_x}+\|u_0\|^{1+\omega}_{H^{2}_x}.
\end{align*}
Thus (\ref{pi6}) and the bound $\|u\|_{H^1_x}\le C(\|u_0\|_{H^1_x})$ yield
\begin{align}\label{pi8}
\|u(t)\|^2_{{H}^2_x}&\le C(\|u_0\|_{H^1_x})(\|u_0\|_{H^{2}_x}+1)
\end{align}
It is standard that if on any interval of the form $t\in[\tau,\tau+T]$ with an fixed $T$ depending only on $\|u_0\|_{H^1_x}$ there holds
\begin{align}\label{pi8}
\|u(t)\|^2_{{H}^2_x} &\le C(\|u_0\|_{H^1_x})(\|u(\tau)\|_{H^{2}_x}+1),
\end{align}
then one has the polynomial bound
\begin{align}\label{pi8}
\|u(t)\|_{{H}^2_x}&\le C(\|u_0\|_{H^1_x})(t+1).
\end{align}
One may see Bourgain \cite{Bo1,Bo2} and Staffilani \cite{Sta} for this.
Therefore, (\ref{ui}) follows due to the time invariance and the global $\|u(t)\|_{H^1}$ norm.
\end{proof}

\begin{Lemma}
Denote $B(s) = s{\left( {K} \right)^{\frac{s}{2}}} + [ {x \cdot \nabla ,{{\left( {K} \right)}^{\frac{s}{2}}}}],$ then
$$B(s)=c\int_{\Gamma}\tau^{\frac{s}{2}}(\tau-K)^{-1}V_3(\tau-K)^{-1}d\tau.$$
\end{Lemma}
\begin{proof}
The proof is adapted from Lemma 7.1 in \cite{CGV}. Recall $K=\mathcal{H}(\alpha_1)$ and the definition of $\mathcal{H}(\alpha_1)^{\frac{s}{2}}$ in Section 2, then
\begin{align*}
2\pi iB(s)f = s\int_\Gamma  {{\lambda ^{\frac{s}{2} - 1}}} {(\lambda  - K)^{ - 1}}Kfd\lambda + \int_\Gamma  {{\lambda ^{\frac{s}{2} - 1}}} \left[ {x{\partial _x},{{(\lambda  - K)}^{ - 1}}K} \right]fd\lambda.
\end{align*}
It is easy to verify
\begin{align*}
\left[ {x{\partial _x},{{(\lambda  - K)}^{ - 1}}K} \right] = -2\lambda K{(\lambda  - K)^{ - 2}} + \lambda {(\lambda  - K)^{ - 1}}{V_3}{(\lambda  - K)^{ - 1}}.
\end{align*}
Then we obtain
\begin{align*}
 2\pi iB(s)f =& s\int_\Gamma  {{\lambda ^{\frac{s}{2} - 1}}} {(\lambda  - K)^{ - 1}}Kfd\lambda  - 2\int_\Gamma  {{\lambda ^{\frac{s}{2}}}} K{(\lambda  - K)^{ - 2}}fd\lambda  \\
 &+ \int_\Gamma  {{\lambda ^{\frac{s}{2}}}{{(\lambda  - K)}^{ - 1}}{V_3}{{(\lambda  - K)}^{ - 1}}fd\lambda }.
\end{align*}
It remains to prove
\begin{align}\label{huhuasy}
s\int_\Gamma  {{\lambda ^{\frac{s}{2} - 1}}} {(\lambda  - K)^{ - 1}}Kfd\lambda - 2\int_\Gamma  {{\lambda ^{\frac{s}{2} - 1}}} \lambda K{(\lambda  - K)^{ - 2}}fd\lambda=0.
\end{align}
By the identity,
$$2{\lambda ^{\frac{s}{2}}}K{(\lambda  - K)^{ - 2}} =  - 2\frac{d}{{d\lambda }}\left( {{\lambda ^{\frac{s}{2}}}K{{(\lambda  - K)}^{ - 1}}} \right) + s{\lambda ^{\frac{s}{2} - 1}}K{(\lambda  - K)^{ - 1}},
$$
and Lemma \ref{41}, one easily obtains (\ref{huhuasy}), thus finishing our proof.
\end{proof}

\end{document}